\documentclass[a4paper,11pt] {article}
\usepackage[latin1]{inputenc}
\usepackage[T1]{fontenc}
\usepackage[english]{babel}
\usepackage{amsmath}
\usepackage{amsfonts}
\usepackage{amssymb}
\usepackage{graphicx}
\usepackage{float}
\usepackage{multicol}
\usepackage{color}
\definecolor{gris25}{gray}{0.55}
\usepackage{shadow}
\usepackage{makeidx}
\newcommand{\be}{\begin{equation}}
\newcommand{\ee}{\end{equation}}
\newcommand{\bd}{\begin{displaymath}}
\newcommand{\ed}{\end{displaymath}}

\newcommand{\ba}{\begin{eqnarray}}
\newcommand{\ea}{\end{eqnarray}}
\newcommand{\ban}{\begin{eqnarray*}}
\newcommand{\ean}{\end{eqnarray*}}

\newcommand{\R} {\mathbb{R}}
\newcommand{\ZZ} {\mathbb{Z}}

\newcommand{\E} {\mathbb{E}}

\newcommand{\I} {\mathbb{I}}
\newcommand{\Z} {{\bf Z}}
\newcommand{\cov} {\mbox{cov}}

\newcommand{\IRS} {\mbox{IRS}}
\newcommand{\ds}{\displaystyle}
\newenvironment{remark}{\textbf{Remark: }}{}
\newenvironment{example}{\textbf{Example: }}{}
\newenvironment{examples}{\textbf{Examples: }}{}

\def\egalelaw{\renewcommand{\arraystretch}{0.5}
\begin{array}[t]{c}\stackrel{{\cal D}}{\sim} \\
{}\end{array}\renewcommand{\arraystretch}{1}}

\def\convlaw{\renewcommand{\arraystretch}{0.5}
\begin{array}[t]{c}\stackrel{{\cal D}}{\rightarrow} \\
{}\end{array}\renewcommand{\arraystretch}{1}}

\newcommand{\Keywords}[1]{\par\noindent
{\small{ \textbf{\emph{Keywords}}\/}: #1}}

\usepackage{theorem}
\theoremstyle{break}
\newtheorem{theo}{Theorem}
\newtheorem{defi}{Definition}
\newtheorem{prop}{Property}
\newtheorem{lemm}{Lemma}

\theorembodyfont{\rmfamily}
\newtheorem{prfT}{Proof of Theorem}

\newtheorem{prfL}{Proof of Lemma}
\newtheorem{prfP}{Proof of Property}

\begin{document}
\title{Local estimation of the Hurst index of multifractional Brownian motion by Increment Ratio Statistic method}

\date{ }
\maketitle \vspace{-1.5cm}
 \begin{center}
 Pierre, R.~BERTRAND${}^{1,2}$, Mehdi FHIMA${}^{2}$ and Arnaud GUILLIN${}^{2}$\\
 ${}^{1}$  {\it INRIA Saclay} \\
${}^{2}$ {\it Laboratoire de Math\'ematiques, UMR CNRS 6620\\
\& Universit\`e de Clermont-Ferrand II, France}
\end{center}

\begin{abstract}

We investigate here the Central Limit Theorem of the Increment
Ratio Statistic of a multifractional Brownian motion, leading to a
CLT for the time varying Hurst index. The proofs are quite simple
relying on Breuer-Major theorems and an original {\it freezing of
time} strategy. A simulation study shows the goodness of fit of this estimator.\\


\end{abstract}

\section*{Introduction}
The aim of this paper is a simple proof of Central Limit Theorem
(CLT in all the sequel) for the convergence of Increment Ratio
Statistic  method (IRS in all the sequel) to a time varying Hurst
index.

Hurst index is the main parameter of fractional Brownian motion
(fBm in all the sequel), it belongs to the interval $(0,1)$ and it
will be denote by $H$ in all the following. For fBm, the Hurst
index drives both path  roughness, self-similarity and long memory
properties of the process. FBm was introduced by Kolmogorov
\cite{Kolmogorov:1940}  as Gaussian "spirals" in Hilbert space
 and then popularized by Mandelbrot \& Van Ness
\cite{Mandelbrot:VanNess:1968}  for its relevance in many
applications. However, during the two last decades, new devices
have allowed access to large then huge datatsets. This put in
light that fBm itself is a theoretical model and that in real life
situation the Hurst index is, at least, time varying. This model,
called multifractional Brownian motion (mBm) has been introduced,
independently by L\'evy-V\'ehel \& Peltier
\cite{LevyVehel:Peltier:1996} and Benassi \emph{et al}
\cite{Benassi:etal:1997}. Other generalizations of fBm remain
possible, for {\it e.g.} Gaussian processes with a Hurst index
depending of the scale, so-called multiscale fBm \cite{Bardet:Bertrand:2007b}, when $H$ is piecewise constant as in the Step
Fractional Brownian Motion see \cite{ayache:etal:2007}, or a wide
range of Gaussian or non-Gaussian processes fitted to applications
(see for example \cite{cheridito:2003, Bardet:Bertrand:2007}).
However, for simplicity of the presentation, in this work we
restrict ourselves to mBm.

In statistical applications, we estimate the time varying Hurst
index through a CLT. Actually, CLT provides us confidence
intervals. Different statistics can be used to estimate the Hurst
index. Among the popular methods, let us mention quadratic
variations, generalized quadratic variations, see
\cite{Benassi:etal:1998,Coeurjolly:2001,Coeurjolly:2005},  and wavelet analysis, see
e.g \cite{Abry:etal:2003} or \cite{Bardet:Bertrand:2010}. Above
methods  can be expansive in term of time complexity. For this
reason, Surgailis \emph{et al} \cite{Surgailis:etal:2008} and
Bardet \& Surgailis \cite{Bardet:Surgailis:2009} have proposed a
new statistic named increment ratio which can be used for
estimating the Hurst index $H$ and  is faster than the wavelet or
the quadratic variations methods, at the price of a slightly
larger variance.

CLT for the different estimators of Hurst index are presently
standard in the case of fBm, but became very technical in the case
of mBm.  The main novelty of our work is the simplicity of the
proofs. In our point of view, mBm is a fBm where the constant
Hurst index $H$ has been replaced by time varying Hurst index.
It is well known that the random field $(H,t)\mapsto B(H,t)$ is
irregular with respect to time $t$, actually with regularity $H$
which belongs to $(0,1)$. It is less kown that this field is
infinitely differentiable with respect to $H$, see
Meyer-Sellan-Taqqu (1999) and Ayache and Taqqu (2005). Thus, for
all time $t_0$, we can freeze the time varying Hurst index, and
the mBm behaves approximatively like a fBm. Eventually, CLT for
mBm follows from CLT for fBm combined with a control of "freezing
error". This new and natural technology allows us to go further
and obtain for example a CLT for the Hurst function evaluated at a
finite collection of times and also quantitative convergence speed
in the CLT. Note that, up to our knowledge, the "freezing Hust
index" strategy for estimation in mBm was introduced, without
further proof, in Bertrand {\it et al} \cite{Bertrand:etal:2010}.
\\

The remainder of this paper is organized as follows. In
Section~\ref{sec1}, we recall a definition of fBm and the
definition of the Increment Ration Statistic. Next, in
section~\ref{sec2}, we review definitions of fBm and mBm and
precise the localization procedure (or freezing). The main result
is stated in Section~\ref{sec3} and some numerical simulations are
presented in Section~\ref{sec4}.
All technical proofs are postponed in Section 5. \\


\section{Recall on fBm and  Increment Ratio Statistic}
\label{sec1} In this section, we present the Increment Ratio
Statistic (IRS) method obtained by filtering centered Gaussian
processes with stationary increments. Before, we recall definition
of the processes under consideration.
\subsection{Definition of fBm and Gaussian processes with stationary increments}
We describe fBm through its harmonizable representation. However,
it is simpler to adopt a more general framework and then specify
fBm as a particular case. Let $X=\left (X(t), t \in [0,1]\right)$
be a zero mean Gaussian process with stationary increments, the
spectral representation theorem (see  Cram\`{e}r and Leadbetter
\cite{Cramer:Leadbetter:1967} or Yaglom \cite{Yaglom:1957}),
asserts that the following representation is in force
 \ban
 \label{repr:accr:stat:harmonizable}
X(t)&=& \int_{\R} (1-e^{it\xi}) \cdot f^{1/2}(\xi) \;
dW(\xi),~~~~\mbox{for all}~~t \in [0,1], \ean where $W(dx)$ is a
Wiener measure with adapted real and imaginary part such that
$X(t)$ is real valued for all $t$. The function $f$ is a Borelian
even, positive and is called spectral density of $X$. To insure
convergence of the stochastic integral, $f$ should satisfies the
condition given by \ba \label{cond:f}
 \int_{R}
 \big(1\wedge|\xi|^2\big)\cdot
f (\xi) \; d \xi < \infty. \ea

\example Fractional Brownian motion with Hurst parameter $H \in (0,1)$ and scale parameter $\sigma>0$
corresponds to a spectral density given by \ba \label{densiteSpec}
f(\xi)=\sigma^2 |\xi|^{-(2H+1)} \text{ for all } \xi \in \R. \ea
In this paper, we denote fBm by $B(H,t)$ when $\sigma=1$. Stress
that this choice is not the conventional one. But,  IRS is
homogeneous and does not depends on a multiplicative factor. Thus,
in sake of simplicity, we can impose the extra condition
$\sigma=1$.
\subsection{Definition of the $a$-Generalized increments}
In all the sequel, we consider the observation of the process $X$
at discrete regularly spaced times, that is the observation of
$\left(X(t_0), \ldots, X(t_n)\right)$
 at times $t_k=k/n$.
Secondly, we consider a filter denoted by $a$ of length $L + 1$ and
of order $p \geq 1$, where $p \leq L$ are two integers. It
corresponds to an arbitrary finite fixed real sequence $a:=(a_0,
\ldots,a_L)\in \R^{L+1}$ having $p$ vanishing moments, i.e., \ba
\label{filter} \sum_{l=0}^{L} a_l l^i = \left\{
\begin{array}{ll}
    \ds 0 & \mbox{if } i \in \{0,\ldots,p-1\} \\
    \ds \sum_{l=0}^{L} a_l l^p  \neq 0 & \mbox{if } i=p.\\
\end{array}
\right. \ea
\\ Consequently, it is easy to prove, for any integer
$m$, that \ba \label{doubleSum} \sum_{l_1=0}^{L}\sum_{l_2=0}^{L}
a_{l_1}a_{l_2}\left| l_1 - l_2 \right|^{m} = \left\{
\begin{array}{ll}
    \ds 0 & \mbox{if } m \in \{0,\ldots,2p-1\} \\
    \ds \left(\sum_{l=0}^{L} a_l l^p \right)^2   \neq 0 & \mbox{if } m = 2p.\\
\end{array}
\right.\ea The family of such filters will be denotes
$\mathcal{A}(p,L)$.
Then, the $a$-Generalized increments of the discrete
process $\left( X(t_k) \right)_{0 \leq k \leq n}$ are defined, for
all $0 \leq k \leq n-L-1$, as follows \ba \label{delta:tk}
\Delta_a X(t_k) = \sum_{k=0}^{L}a_l X(t_{k+l}) \ea and their
harmonizable representations are given by \ban \Delta_a X(t_k) =
\int_{\R} e^{i t_k \xi} g_a(-\xi/n)f^{1/2}(\xi) \; dW(\xi)
 \ean where $g_a(\cdot)$ is specified as follows \ba \label{g_a} g_a(u):=\sum_{l=0}^{L}a_{l}e^{i l u}.\ea
\text{}\\ \examples In the simple case where $a:=(a_0=1,a_1=-1)$,
the operator $\Delta_a$ corresponds to a discrete increment of
order 1, and when $a:=(a_0=1,a_1=-2,a_2=1)$, the operator
$\Delta_a$ represents the second order differences.

\subsection{Definition of the Increment Ratio Statistic}
%
Let $\left ( \Delta_a X(t_k) \right )_{0 \leq k \leq n-L-1}$ be
the $a$-Generalized increments sequence defined by
$\eqref{delta:tk}$ from the discrete observation
$\left( X(t_k) \right)_{0
\leq k \leq n}$. Then, the IRS introduced by Bardet and Surgailis
\cite{Bardet:Surgailis:2009} is given by

\ba \IRS_{a,n}(X)= \frac{1}{n-L} \sum_{k=0}^{n-L-1} \psi \left(
\Delta_a X(t_k) , \Delta_a X(t_{k+1}) \right) \label{IRSgaus} \ea
where $\psi(\cdot,\cdot)$ is described as follows

\ban \psi(x,y) := \left\{
\begin{array}{ll}
    \ds \frac{\left|x+y \right|}{\left| x \right|+\left|y \right|}& \mbox{if } (x,y) \in \R^2\backslash \{(0,0)\} \\
    \ds \       1 & \mbox{if } (x,y)=(0,0).\\
\end{array}
\right. \ean
\paragraph{IRS of fractional Brownian motion}
\text{}\\
\text{}\\
In the case of the fBm with Hurst parameter $H \in \left (0,1
\right)$, \emph{i.e} $X(t)=B_H(t)$, Bardet and Surgailis have
established in \cite[Corollary 4.3, p.13]{Bardet:Surgailis:2009},
under some semi-parametric assumptions, the following CLT for the
statistics $\IRS_{a,n}$ \ba \label{IRSfBm}
    \sqrt{n} \left( \IRS_{a,n}(B_H)-\Lambda_a(H) \right) \convlaw
    \mathcal{N}(0,\Sigma^2_{a}(H)) \text{ with } \left\{
\begin{array}{ll}
    \ H \in \left (0,3/4 \right) & \mbox{if } a=(1,-1)\\
    \ H \in \left (0, 1 \right) & \mbox{if } a=(1,-2,1)\\
\end{array}
\right. \label{tclB}\ea where the sign $\convlaw$ means
convergence in distribution,
 \ba  \label{rhoH}
 \Lambda_a(H) &:=& \Lambda_0\left( \rho_a (H) \right) \label{Lp}  \\
 \Lambda_0(r) &:=&  \frac{1}{\pi} \arccos (-r)
 + \frac{1}{\pi} \sqrt{\frac{1+r}{1-r}} \log \left( \frac{2}{1+r} \right) \label{L0} \\
 \rho_a (H) &=& \left\{
\begin{array}{ll}
    \ds 2^{2H-1}-1 & \mbox{if } a=(1,-1)\\
    \ds \frac{-3^{2H}+2^{2H+2}-7}{8-2^{2H+1}} & \mbox{if } a=(1,-2,1)\\
\end{array}
\right. \ea and the asymptotic variance $\Sigma^2_{a}(H)$ is given
by \ban \Sigma^2_{a}(H)=\sum_{j \in \ZZ} \cov \left( \psi\left(
\Delta_a B_H(t_0) , \Delta_a B_H(t_{1}) \right), \psi \left(
\Delta_a B_H(t_j) , \Delta_a B_H(t_{j+1}) \right) \right)
\label{Sigma_a}.
 \ean The graphs of $\Lambda_0(\rho), \rho_a(H)$ and $\Lambda_a(H)$, with $a=(1,-1) \text{ or } a=(1,-2,1)$,
 are given in Figure~\ref{fig1}, Figure~\ref{fig2} and Figure~\ref{fig3}
 below.
\begin{figure}[htbp]
\begin{center}
\begin{tabular}{c}
\includegraphics[width=8cm,height=3cm]{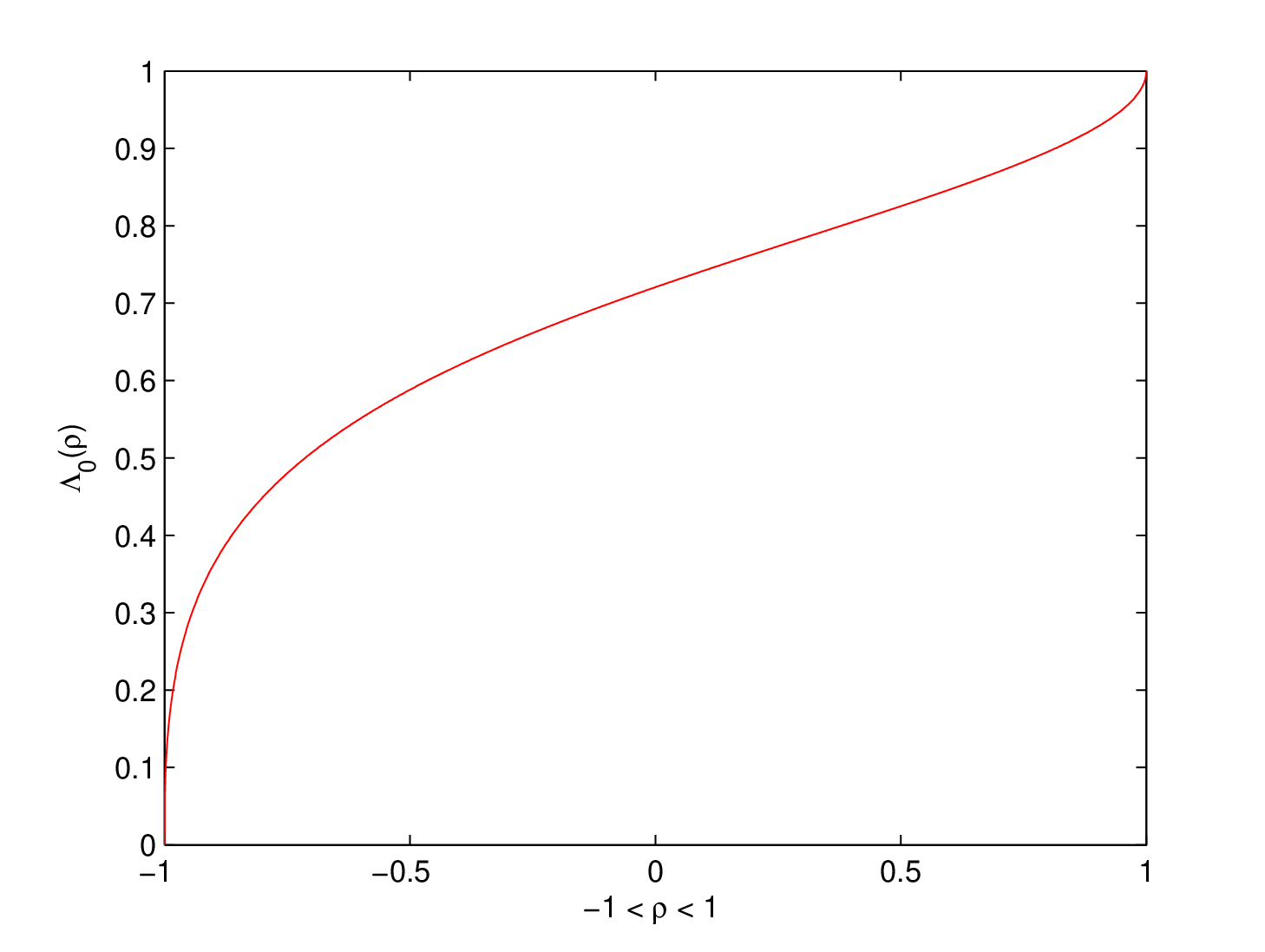}
\end{tabular}
\caption{\emph{The graph of $\Lambda_0(\rho)$.}} \label{fig1}
\end{center}
\end{figure}
\begin{figure}[htbp]
\begin{center}
\begin{tabular}{c}
\includegraphics[width=8cm,height=3cm]{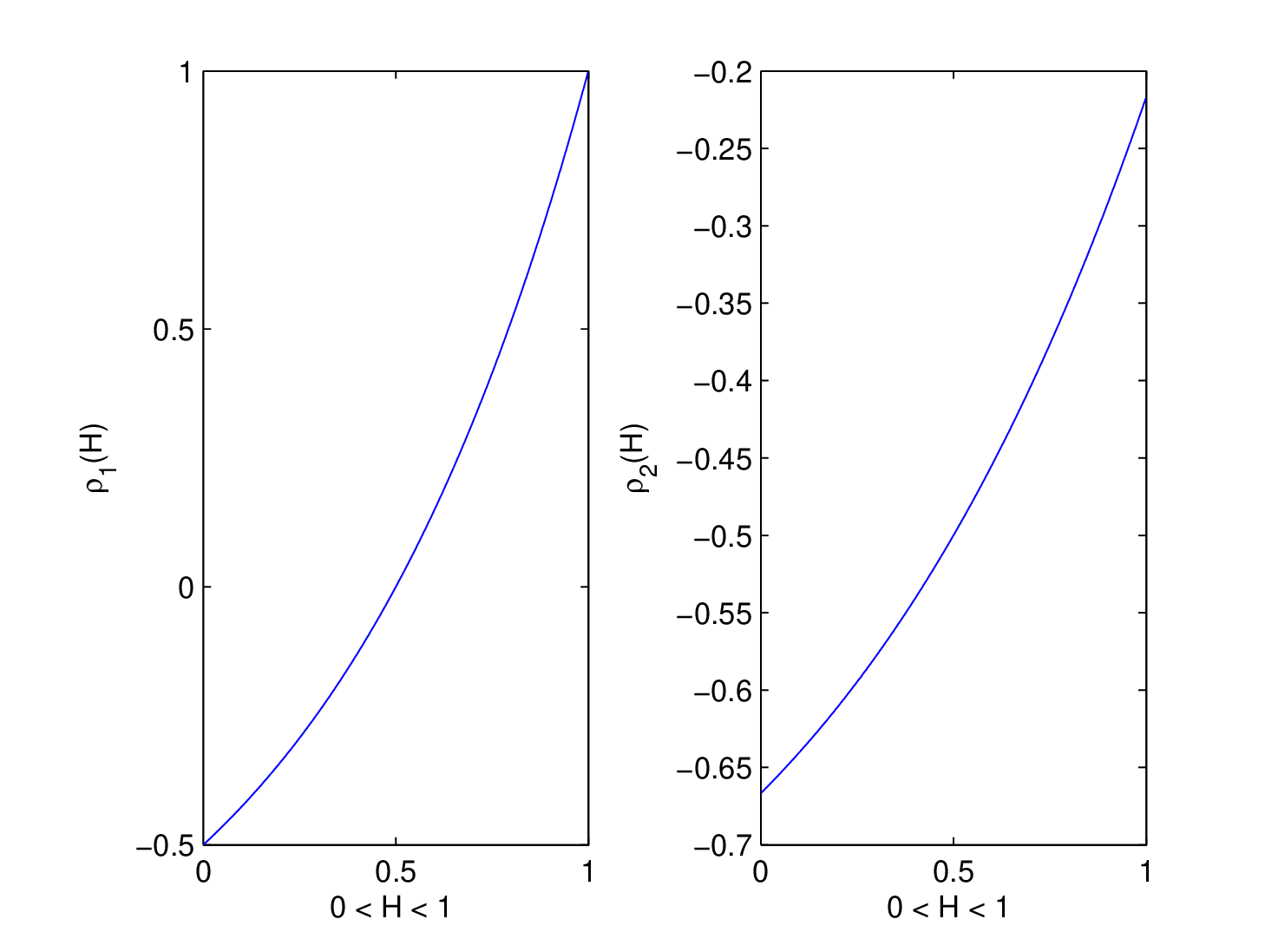}
\end{tabular}
\caption{\emph{The graphs of $\rho_a(H)$ with $a=(-1,1)$(left) and
$a=(1,-2,1)$ (right).}} \label{fig2}
\end{center}
\end{figure}
\begin{figure}[htbp]
\begin{center}
\begin{tabular}{c}
\includegraphics[width=8cm,height=3cm]{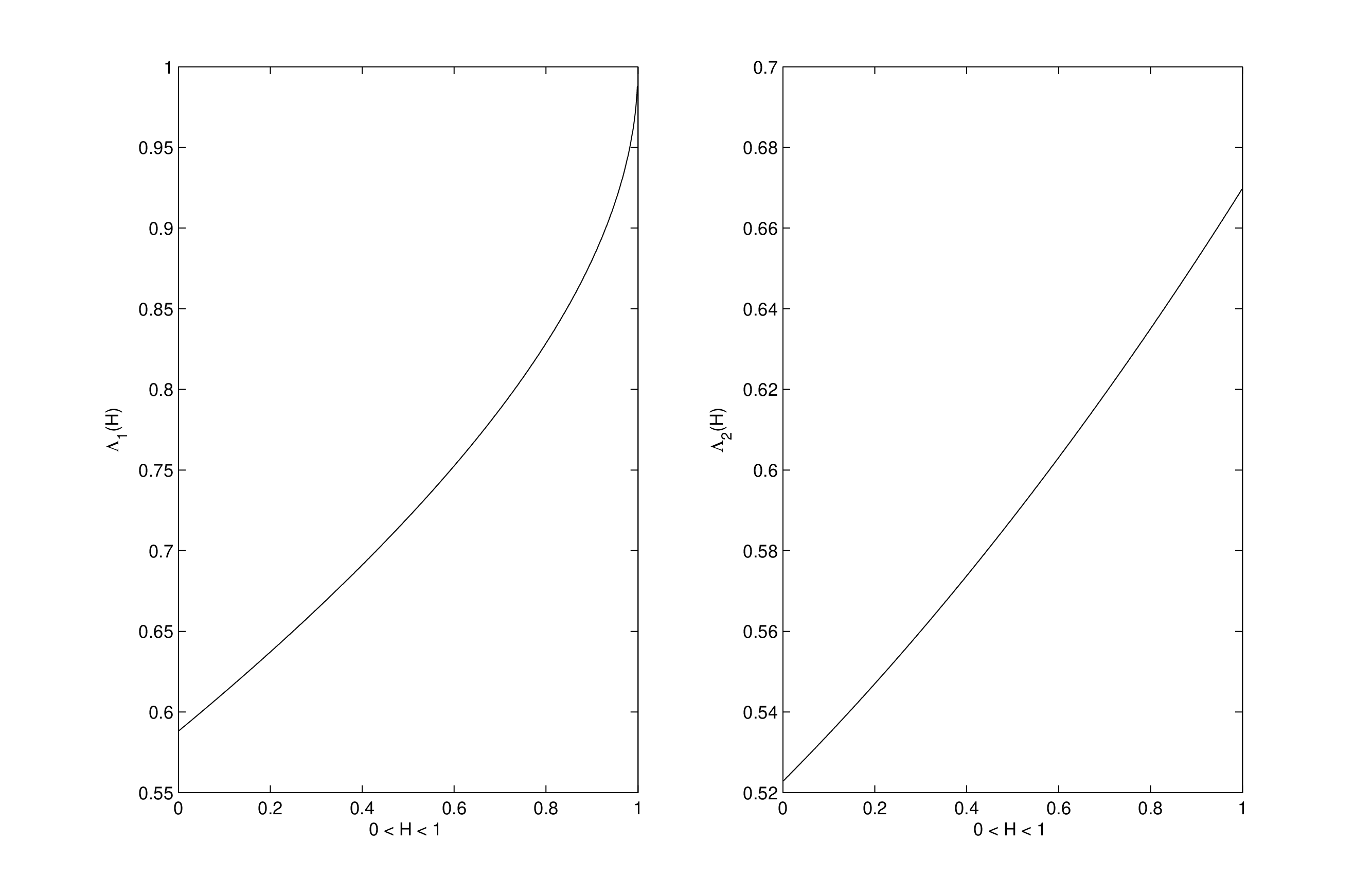}
\end{tabular}
\caption{\emph{The graphs of $\Lambda_a(H)$ with $a=(-1,1)$(left)
and $a=(1,-2,1)$  (right).}} \label{fig3}
\end{center}
\end{figure}
It is easy to prove that the function $H \mapsto
\Lambda_a(H)$, with $a=(1,-2,1)$, is a monotonic increasing
function  in the interval (0,1), see Figure~\ref{fig3}. Therefore,
$\ds \widehat{H}_n=\Lambda^{-1}_a \left ( \IRS_{a,n}(B_H)\right)$
provides  an estimator of the Hurst parameter $H$ with convergence rate $
\mathcal{O }\left( \sqrt{n} \right)$.
 Moreover, we refer to
Stoncelis and Vai\v{c}iulis \cite{Stoncelis:Vaiciulis:2008} for a
numerical approximation of the variance $\Sigma^2_a(H)$ with
$a=(1,-1) \text{ or } a=(1,-2,1)$, needed for construction of
confidence intervals, see \cite[Corollary 4.3,
p.13 and Appendix, p.32]{Bardet:Surgailis:2009}.

\section{Going from fBm to mBm and return by  freezing}\label{sec2}
The main goal of this section is to present different
representations for the FBm and the mBm enabling us to present the
{\it time freezing} strategy we will use to prove our main
theorems.

\subsection*{FBm and its different representations}
Fractional Brownian motion was introduced by Kolmogorov
\cite{Kolmogorov:1940} and then made popular by Mandelbrot \& Van
Ness \cite{Mandelbrot:VanNess:1968}. This process has been widely
used in applications to model data that exhibit self-similarity,
stationarity of increments, and long range dependence. FBm with
Hurst parameter $H \in(0, 1)$, denoted by $\left( B_H(t), t \in
[0,1] \right)$, is a centered Gaussian process with covariance
function defined for $s, t \in [0, 1]$ by \ba \label{covFBM} \E
\left[B_H(t)B_H(s) \right] = \frac 1 2 \left( t^{2H} + s^{2H} +
|t-s|^{2H} \right ).\ea This process is characterized by its Hurst
index which drives both pathwise regularity, self-similarity and
long memory, see {\it e.g.} the overview in Bertrand et al (2010).
Before going further, let us precise notations: in all the sequel
we will denote by $ B_H$ the fBm and  $B(H,t)$ the random field
depending on both Hurst index and time.  Up to a multiplicative
constant the two notions coincide, more precisely we have $ B_H(t)
= C(H) \times B(H,t)$ for $C(H)$ a non-negative constant depending
on $H$.

Fractional Brownian motion, $\left( B_H(t), t \in [0,1] \right)$,
can be represented through its harmonizable representation
(\ref{repr:accr:stat:harmonizable}, \ref{densiteSpec}), or its
moving-average representation (see Samorodnitsky \& Taqqu
\cite[Chapter 14]{Samorodnitsky:Taqqu:1994}). A third
representation is the wavelet series expansion introduced by Meyer
\emph{et al} \cite{Meyer:etal:1999}, and then nicely used by
Ayache and Taqqu (2003 and 2005). In this case, it is judicious
to shift to the random field representation defined as follows \ba
\label{waveletSeries} B(H,t)&=&\sum_{j \in Z} \sum_{k \in \Z}
a_{jk}(t,H)\,\epsilon_{jk},~~~~\mbox{for all}~~t \in [0,1]
 \ea
 where
 $(\epsilon_{jk})_{(j,k) \in Z^2}$ is a sequence of standard
Gaussian random variables $\mathcal{N}(0,1)$, the non-random
coefficients $a_{jk}(t,H)$ are given by $\ds a_{jk}(t,H)=\int_{\R}
(1-e^{it\xi}) \cdot |\xi|^{-(H+1/2)} \widehat{\psi}_{j,k}(\xi) \;
d\xi$, and $\widehat{\psi}_{j,k}$ is the Fourier transform of the
Lemarié-Meyer wavelet basis $\psi_{j,k}$. Let us refer to Ayache
and Taqqu (2003) for all the technical details. To put it into a
nutshell, by using the Meyer \emph{et al} 's Lemma
(\cite{Meyer:etal:1999}), we can prove the existence of an almost
sure event $\Omega^*$, that is such that $\Pr(\Omega^*)=1$, such
that for all $\omega \in \Omega^*$ the series (\ref{waveletSeries}
) converges uniformly for $(H,t)\in K$ where $K$ is any compact
subset of $(0,1)\times \R$. Moreover, the field is infinitively
differentiable with respect to $H$ with derivatives bounded
uniformly on every compact subset of $(0,1)\times \R$ by a
constant $C^*(\omega)>0$ where $C^*$ is  a positive random
variable with finite moments of every order.

\subsection*{MBm and its different representations}
To be short, mBm is obtained by plugging a time varying Hurst
index $t\mapsto H(t)$ into one of the three representations of the
fBm given above, that is the moving average representation, the
harmonizable one ((\ref{repr:accr:stat:harmonizable},
\ref{densiteSpec})) or the wavelet series expansion
(\ref{waveletSeries} ). The function $t\mapsto H(t)$ should be at
least continuous, and if the Hölder regularity of the function
$t\mapsto H(t)$ is greater than $\max(H(t), t\in [0, 1])$  (the
so-called condition $\bf(\mathcal{C})$ in Ayache and Taqqu
(2003)), then for every time $t\in[0, 1]$ the roughness of mBm is
given by $H(t)$. Les us also refer to Cohen \cite{Cohen:1999}
where he proves that the moving average representation and the
harmonizable representation of mBm are equivalent up to a
multiplicative deterministic function, and to Meyer {\it et al.}
to the almost sure equality of harmonizable representation and
wavelet series expansion.
\par
\medskip

With this tools, we are now in order to precise our "freezing" technology:

\subsection*{MBm behaves locally as a fBm}
By applying Taylor expansion of order 1 around any fixed time $t^*
\in [0,1]$, we obtain the following formula \ba \label{taylor}
B(H(t),t) \I_{\Omega^{*}}&=& B(H^*,t) \I_{\Omega^{*}} + R(t)
\I_{\Omega^{*}} \ea where $R(t)$ refers to the Taylor rest which
satisfies \ba \label{rest} \sup_{s \in [0,1]} \left |R(s)
\I_{\Omega^{*}} \right | \leq C^{*}(\omega) |H(t)-H^{*}| \ea with
$C^{*} >0$ a positive random variable with finite moments of every
order. Noting that $H^{*}$ corresponds to the value of the Hurst
function $H(\cdot)$ at $t^*$, i.e., $H^*=H(t^*)$ and
$\I_{\Omega^{*}}$ represents the indicator function of a subset
$\Omega^{*}$ defined by: $\I_{\Omega^{*}}(\omega) = \left\{
\begin{array}{ll}
    \ 1 & \mbox{if } \omega \in \Omega^{*} \\
    \ 0 & \mbox{else }
\end{array}
\right.$. Next, if we know that the Hurst function $H(\cdot)$ has a
Hölder regularity of order $\eta >0$, so we obtain immediately
that \ba \label{supR} \sup_{s \in [0,1]} \left |R(s)
\I_{\Omega^{*}} \right | \leq M^{*}\, \I_{\Omega^{*}}\,
|t-t^{*}|^{\eta} \ea with $ M^{*}=c \times C^* > 0$ a positive random
variable with finite moments of every order.
\section{Main results}
\label{sec3}

This section is dedicated to the CLT of the IRS localized version
for the mBm. Let us however first give a simple result on the CLT
for the IRS of Gaussian processes with stationary increments,
which, applied to the fractional Brownian motion, gives with a
simple proof the result of Bardet-Surgailis
\cite{Bardet:Surgailis:2009}.

We thus consider a process $X$ observed through the knowledge of
$\left(X(t_0), \ldots, X(t_n)\right)$ with $t_k=k/n$ for $k=1,
\dots, n$. The corresponding increment ratio statistic
$\IRS_{a,n}(X)$ is defined by $\eqref{IRSgaus}$, with a filter  $a
\in \mathcal{A}(p,L)$, that is satisfying $\eqref{filter}$.
\begin{theo}[Fractional Brownian motion]
\label{th:Theorem1}
\begin{itemize}
\item[i)]
Let $X$ be a zero mean Gaussian process with stationary increments.
We assume that \ba \label{condG}
\sum_{j \in \ZZ} |r_{a}(j)|^2 <+\infty \ea where $r_{a}(j):=\cov
\left( \Delta_a X(t_0), \Delta_a X(t_j) \right)$ for $j \in \ZZ$
is supposed independent of $n$. Then
 \ba \label{TCL:fBm}
 \sqrt{n} \left(
\IRS_{a,n}(X)-\Lambda_0(\rho_{a}) \right)& \convlaw&
\mathcal{N}(0,\Sigma^2_a) \label{TCLGauss} \ea where
$\Lambda_0(\cdot)$ is defined by $\eqref{L0}$, $\rho_{a}$
represents the correlation between two successive $a$-Generalized
increments, and the asymptotic variance $\Sigma^2_a$ is given by
\ban \Sigma^2_a:=\sum_{j \in \ZZ} \cov \left( \psi \left( \Delta_a
X(t_0) , \Delta_a X(t_{1}) \right), \psi \left( \Delta_a X(t_j) ,
\Delta_a X(t_{j+1}) \right) \right)\ean and is well defined and
belongs to $[0,+\infty)$. \item[ii)] In particular, let $X$ be a
fBm, that is $X=B_H$ with Hurst parameter $H \in \left (0,1
\right)$. Moreover, in the case $p=1$ assume the extra assumption
$\ H \in \left (0,3/4 \right)$. Then CLT (\ref{TCL:fBm}) is in
force  where $\Lambda_a(\cdot)=\Lambda_0(\rho_a(\cdot))$ is a
monotonic increasing function of $H$ with $\Lambda_0(\cdot)$,
resp. $\rho_a(\cdot)$ described by $\eqref{L0}$, resp.
$\eqref{rho_aH}$, and the asymptotic variance $\Sigma^2_a$ is
given by \ban \Sigma^2_a:=\sum_{j \in \ZZ} \cov \left( \psi \left(
\Delta_a B_H(t_0) , \Delta_a B_H(t_{1}) \right), \psi \left(
\Delta_a B_H(t_j) , \Delta_a B_H(t_{j+1}) \right) \right)\ean
which   is well defined and belongs to
$[0,+\infty)$.\hfill$\square$
\end{itemize}
\end{theo}
\begin{remark}
\begin{enumerate}
\item In the sequel, in order to could inverse function
$\Lambda_a(\cdot)$, we will assume that the filter $a \in
\mathcal{A}(p,L)$ satisfy $L=p \geq 1$ and $a_l=(-1)^{p-l} \binom
{p} {l}$ for all $l \in \{0,\ldots,p\}$. The class of such filters
will be denoted $\mathcal{B}(p)$ and named binomial filters. This
restriction is motivated by the fact that in the particular case
where $a \in \mathcal{B}(p)$, the correlation function
$\rho_a(\cdot)$ defined by $\eqref{rho_aH}$, is a monotonic
increasing function of $H$, instead of in the general case where
$a \in \mathcal{A}(p,L)$ it is not always true.\item The
regularity of $\Lambda_a(H)$ enables us then to get via the well
known Delta-method the CLT for the Hurst parameter. However no
closed formulae for $\Lambda_a(H)^{-1}$ is available so that the
limiting covariance will be no further explicit. \item We stress
once again that the proof of the theorem is quite simple. Note
also that using recent results of Nourdin {\it et al} \cite[Th.
2.2]{Nourdin:etal:2010}, we even have that there exists a sequence
$\gamma(n)$ decaying to zero such that for all $h\in C^2$ and
$N\sim\mathcal{N}(0,\Sigma^2_a) $
$$\left|\mathbb{E}\left[h\left( \sqrt{n} \left(
\IRS_{a,n}(X)-\Lambda_0(\rho_{a})
\right)\right)-h(N)\right]\right|\le \| h''\|_\infty\,\gamma(n).$$
The precise estimation of $\gamma(n)$ is however out of the scope
of the present paper and will be found in \cite{Fhima:2011}. Using
\cite[Cor. 2.4]{Nourdin:etal:2010}, we also have that the previous
CLT may be reinforced to a convergence in 1-Wasserstein distance
or in Kolmogorov distance. \item The reader will have noticed that
the assumption  $r_{a}(j):=\cov \left( \Delta_a X(t_0), \Delta_a
X(t_j) \right)$   independent of $n$ for all $j \in \ZZ$  is a
quite strong one. Indeed, for the multiscale Brownian motion, this
not true. However, in a sense, it is asymptotically true and it
may then be applied to prove the convergence of the IRS to the
Hurst parameter related to the highest frequency. See
\cite{Fhima:2011} for further details.
\end{enumerate}
\end{remark}

\medskip

To achieve our final goal, we state by presenting a Lemma where
we prove that the localized version of IRS for mBm converges in
$L^2(\Omega)$ to the IRS of fBm with a certain rate.\\

\noindent{\bf Localized version of the IRS for multifractional Brownian motion}


Let us consider a multifractional Brownian motion with Hurst
function $H(\cdot)$ denoted by $X=\left(B(H(t),t), t \in
[0,1]\right)$. Secondly, let $t^{*} \in (0,1)$ be an arbitrary
fixed point, then we  denote by $\nu_{n}\left(\gamma,t^{*}\right)$
the set of indices around $t^{*} $, given by \ba \label{voisinage}
\nu_{n}\left(\gamma,t^{*}\right) &=& \{ k \in \{ 0,\ldots,n-L-1 \}
: |t_k-t^{*}| \leq n^{-\gamma} \} \\ &=& \{ \lfloor
nt^{*}-n^{1-\gamma} \rfloor,\ldots, \lfloor nt^{*}+n^{1-\gamma}
\rfloor \} \ea where $\lfloor x \rfloor$ is the integer part of
$x$ and $\gamma \in (0,1)$ is a fixed parameter which allows to
control the size of $\nu_{n}\left(\gamma,t^{*}\right)$ which
cardinal is equal to $v_{n}(\gamma) := 2n^{1-\gamma} +1$. Finally,
for any $n$ large enough, we denote by
$\IRS_{a,n}^{\gamma,t^{*}}\left(B_{H(\cdot)}\right)$ the localized
version of IRS defined as follows \ba \label{IRSmBm}
\IRS_{a,n}^{\gamma,t^{*}}\left(B_{H(\cdot)}\right)&= &\frac{1}{
2n^{1-\gamma} +1} \sum_{k=\lfloor nt^{*}-n^{1-\gamma}
\rfloor}^{\lfloor nt^{*}+n^{1-\gamma} \rfloor} \psi \left (
\Delta_a B_{H(t_k)}(t_k), \Delta_a B_{H(t_{k+1})}(t_{k+1})
\right). \ea With these notations, we are in order to state our
main result:
\begin{theo}[Multifractional Brownian motion]\label{th:Theorem2}
\begin{itemize}
\item[i)]Let $X$ be a mBm, its localized IRS be defined by (\ref{IRSmBm}) and  assume that $\gamma (1+\eta)
>1$. Then \ba \label{TCLmBm} n^{(1-\gamma)/2} \left(
\IRS_{a,n}^{\gamma,t^{*}}\left(B_{H(\cdot)}\right)-\Lambda_a(H^*)
\right) \convlaw \mathcal{N}(0,\Sigma^2_a) \text{ with } \left\{
\begin{array}{ll}
    \ H^* \in \left (0,3/4 \right) & \mbox{if } p=1 \\
    \ H^* \in \left (0, 1 \right) & \mbox{if } p \geq 2 \\
\end{array}
\right.\ea where $\Lambda_a(\cdot)=\Lambda_0(\rho_a(\cdot))$ is a
monotonic increasing function of $H^*$ with $\Lambda_0(\cdot)$ \&
$\rho_a(\cdot)$ described by $\eqref{L0}$ \& $\eqref{rho_aH}$, and
the asymptotic variance $\Sigma^2_a$ is given by \ban
\Sigma^2_a:=\sum_{j \in \ZZ} \cov \left( \psi \left( \Delta_a
B_{H^*}(t_0) , \Delta_a B_{H^*}(t_{1}) \right), \psi \left(
\Delta_a B_{H^*}(t_j) , \Delta_a B_{H^*}(t_{j+1}) \right)
\right)\ean and is well defined and belongs to $[0,+\infty)$.
\item[ii)] Let now consider $0<t_0<t_1<...<t_m$ for a finite $m$
then under the same assumption we can enhance the previous CLT to
the vector
$$n^{(1-\gamma)/2} \left(\IRS_{a,n}^{\gamma,t_1}\left(B_{H(\cdot)}\right)-\Lambda_a(H(t_1)),...,\IRS_{a,n}^{\gamma,t_m}\left(B_{H(\cdot)}\right)-\Lambda_a(H(t_m))\right)$$
with a well defined limiting covariance $\mathfrak{S}$.
\hfill$\square$
\end{itemize}
\end{theo}

\begin{remark}
\begin{enumerate}
\item Here again, one can use results of \cite{Nourdin:etal:2010}
to get explicit estimates on the speed of convergence for this
CLT. \item It is highly interesting to upgrade the previous CLT to
the trajectory level, needing then a tightness result, for example
to test if the Hurst coefficient is always greater than 1/2, or to
perform other test. Such kind of result will be developed in
\cite{Fhima:2011}.
\end{enumerate}
\end{remark}

\section{Numerical results}
\label{sec4} In this section, for numerical estimation of the
Hurst index by IRS, one has chosen a binomial filter of order 2,
i.e. $a=(1,-2,1)$, insuring the convergence of the estimator
$\widehat{H}_n$ for any $H \in (0,1)$. At first, we analyze
through Monte-Carlo simulations the efficiency of the Hurst
parameter of fBm estimator given by IRS. Then, we study the
estimators of some Hurst functions of mBm obtained by localized
version of IRS, and we compare it with the estimators given by
Generalized Quadratic Variations (GQV) method, see e.g Coeurjolly
\cite{Coeurjolly:2005}.

\subsection*{Estimation of the Hurst index of fBm}

At first, by using Wood and Chan \cite{Chen:Wood:1998} algorithm,
for $n=10 000$ we have simulated three replications of the fBm
sequences $B_H=\left(B_H(t_0), \ldots, B_H(t_n)\right)$, at
regularly spaced times such that $t_k=k/n$ with $k=0,\ldots,n$ ,
for three values of the Hurst parameter $H$, denoted
$\{H_1,H_2,H_3 \}$, and given by
\begin{description}
    \item [$\mathbf{(\mathcal{C}_1)}$] $H_1=0.3 < 1/2$ for short
    range dependent case,
    \item [$\mathbf{(\mathcal{C}_2)}$] $H_2=1/2$ for standard
    Brownian motion,
    \item [$\mathbf{(\mathcal{C}_3)}$]$H_3=0.7 > 1/2$ for long
    range dependent case,
\end{description}
see Figure~\ref{fig4} below.\\
\begin{figure}[htbp]
\begin{center}
\begin{tabular}{c}
\includegraphics[width=15cm,height=5cm]{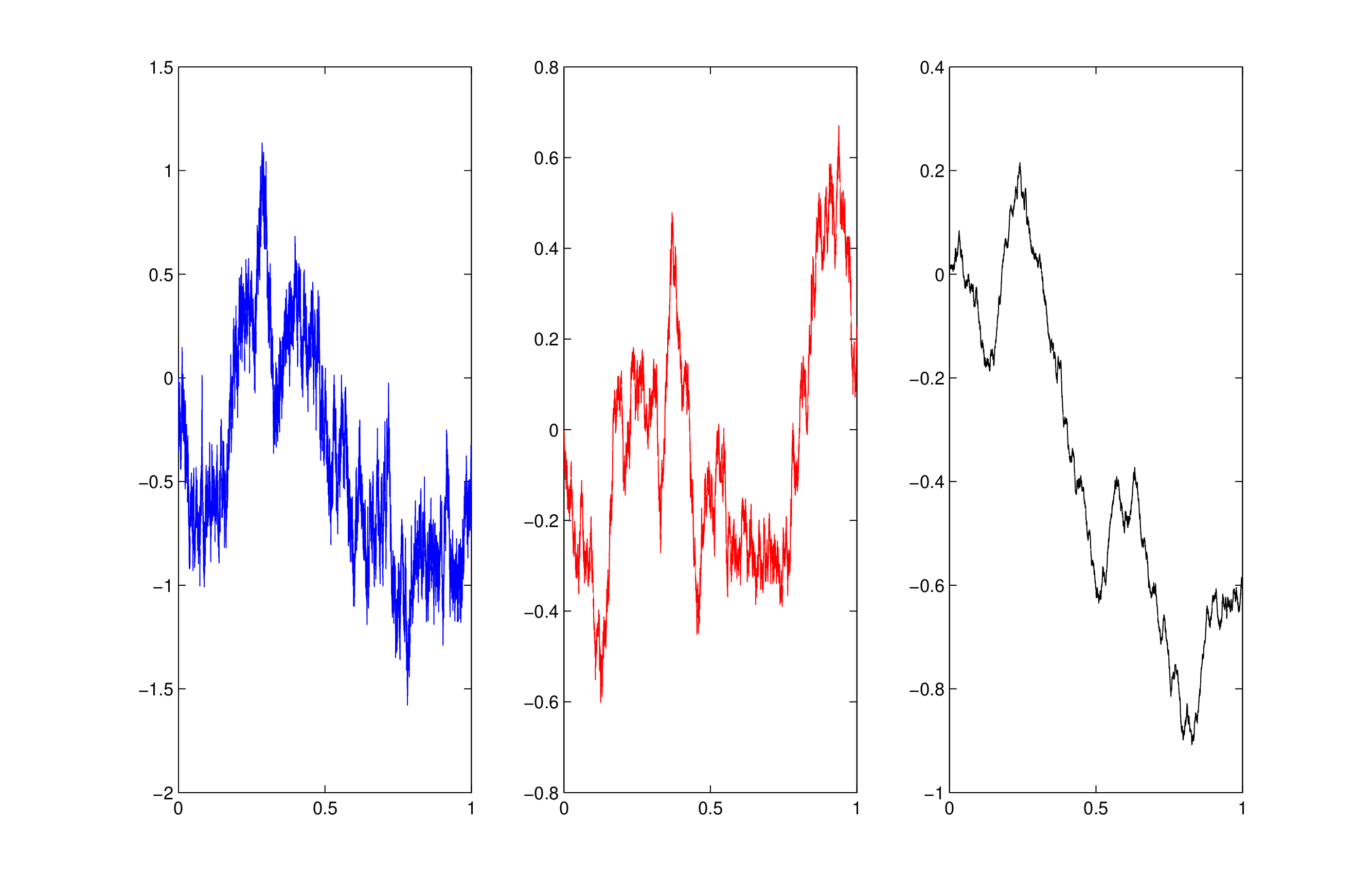}
\end{tabular}
\caption{\emph{Simulated fBm process with $H=0.3$ (left), $H=0.5$
(middle) and $H=0.7$ (right).}} \label{fig4}
\end{center}
\end{figure}
Then, for each sample $\mathbf{(\mathcal{C}_i)}$ with $i \in
[1,3]$, we have computed the increment ratio statistic
$IRS_{2,n}(B_{H_i})$ and estimated the Hurst index given by $\ds
\widehat{H}_{n,i}=\Lambda^{-1}_2 \left (
\IRS_{2,n}(B_{H_i})\right)$. We remark that the IRS methods
provide good results given in Table~\ref{tab1} below.\\
\begin{table}[htbp]
\begin{center}
\begin{tabular}{|c|c|c|c|}
  \hline
  Exact values of $H$ & 0.3 & 0.5 & 0.7 \\
  \hline
  Estimated values of $H$ & 0.3009  & 0.4993 & 0.7000 \\
  \hline
\end{tabular}
\caption{\emph{Estimated values of $H$.}}\label{tab1}
\end{center}
\end{table}
These examples are plainly confirmed by Monte Carlo simulations.
Indeed, for each case $\mathbf{(\mathcal{C}_i)}$ with $i \in
[1,3]$, we have made $M=1000$ simulations of independent copies of
fBm sequences $B^{(k)}_{H_i}=\left(B^{(k)}_{H_i}(t_0), \ldots,
B^{(k)}_{H_i}(t_n)\right)$, for $k=1,\ldots,M$. We find also good
results illustrated by the following histograms, see
Figure~\ref{fig5}, which represent the distribution of the
estimator $\widehat{H}_{n,i}$, for  $i \in [1,3]$. Thus, we have
computed the standard deviation  $\E |\widehat{H} - H|^2$ given in
Table~\ref{tab2} below.\\

\begin{figure}[htbp]
\begin{center}
\begin{tabular}{c}
\includegraphics[width=15cm,height=5cm]{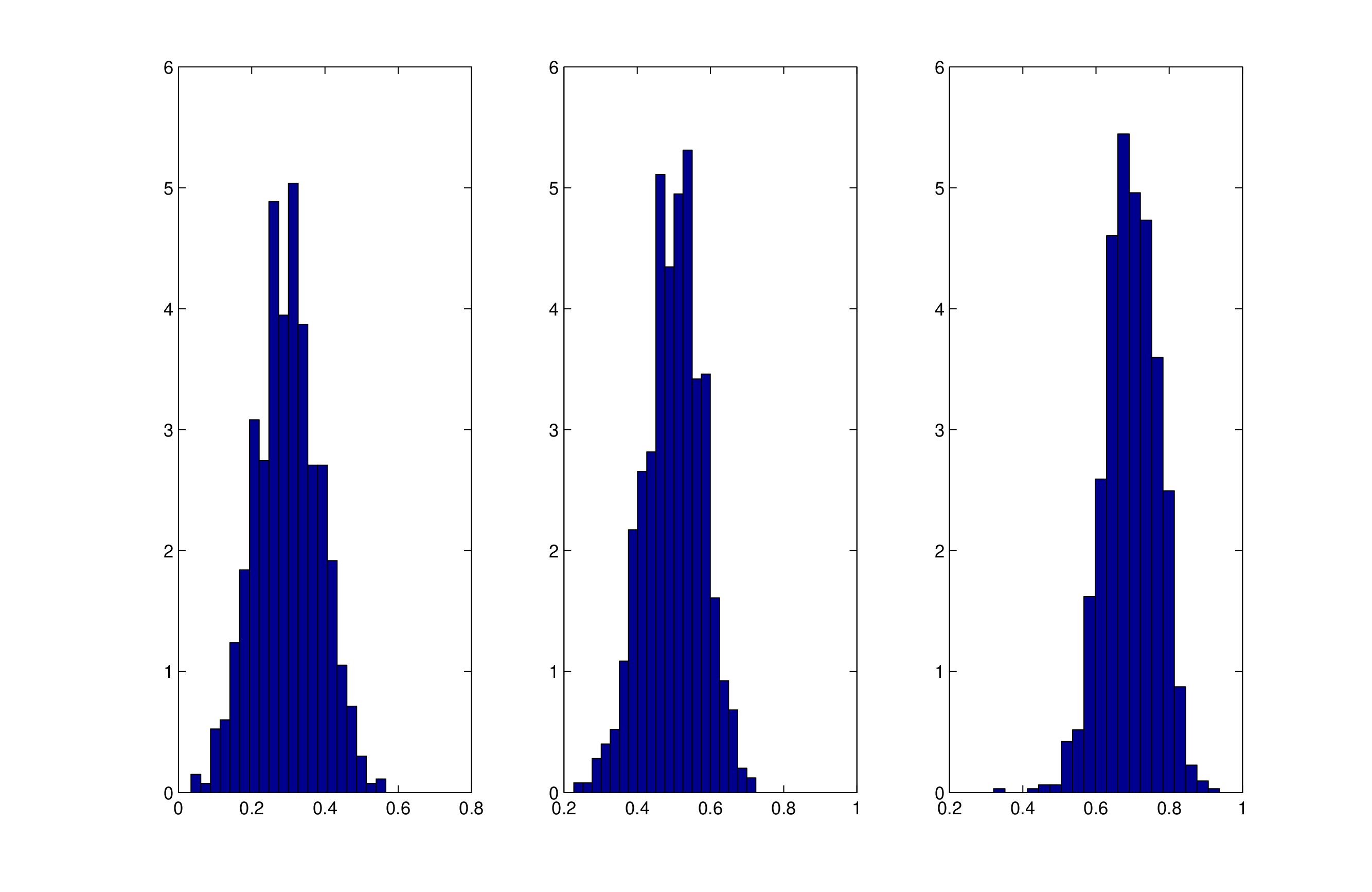}
\end{tabular}
\caption{\emph{Distribution of the estimated values of $H$ in the
case $H=0.3$ (left), $H=0.5$ (middle) and $H=0.7$ (right).}}
\label{fig5}
\end{center}
\end{figure}

\begin{table}[htbp]
\begin{center}
\begin{tabular}{|c|c|c|c|}
  \hline
  Different values of $H$ & 0.3 & 0.5 & 0.7 \\
  \hline
  Standard deviation  & $8.1865 \times 10^{-5}$  & $7.6914 \times 10^{-5}$ & $6.9837 \times 10^{-5}$ \\
  \hline
\end{tabular}
\caption{\emph{Standard deviation.}}\label{tab2}
\end{center}
\end{table}

\subsection*{Local estimation of the Hurst function of mBm}
To synthesis a sample path of a mBm, one has used the Wood and
Chan circulant matrix improved with kriging interpolation method,
which is  faster than Cholesky-Levinson factorization algorithm.
In fact, both methods are not exact but provide good results. For,
$ n = 10 000$, we have simulated three samples of the mBm
sequences $B_{H(\cdot)}=\left(B_{H(t_0)}(t_0), \ldots,
B_{H(t_n)}(t_n)\right)$, at regularly spaced times such that
$t_k=k/n$ with $k=0,\ldots,n$ , for three types of the Hurst
function $H(\cdot)$, namely
\begin{description}
    \item [$\mathbf{(\mathcal{C}_4)}$] Linear function: $H_4(t)=0.1+0.8t$,
    \item [$\mathbf{(\mathcal{C}_5)}$] Periodic function $H_5(t)=0.5+0.3 \sin(\pi t)$,
    \item [$\mathbf{(\mathcal{C}_6)}$] Logistic function: $\ds H_6(t)=0.3+\frac{0.3}{(1+\exp(-100(t-0.7)))}$,
\end{description}
see Figure~\ref{fig6} below.\\
\text{}\\
\begin{figure}[htbp]
\begin{center}
\begin{tabular}{c}
\includegraphics[width=15cm,height=5cm]{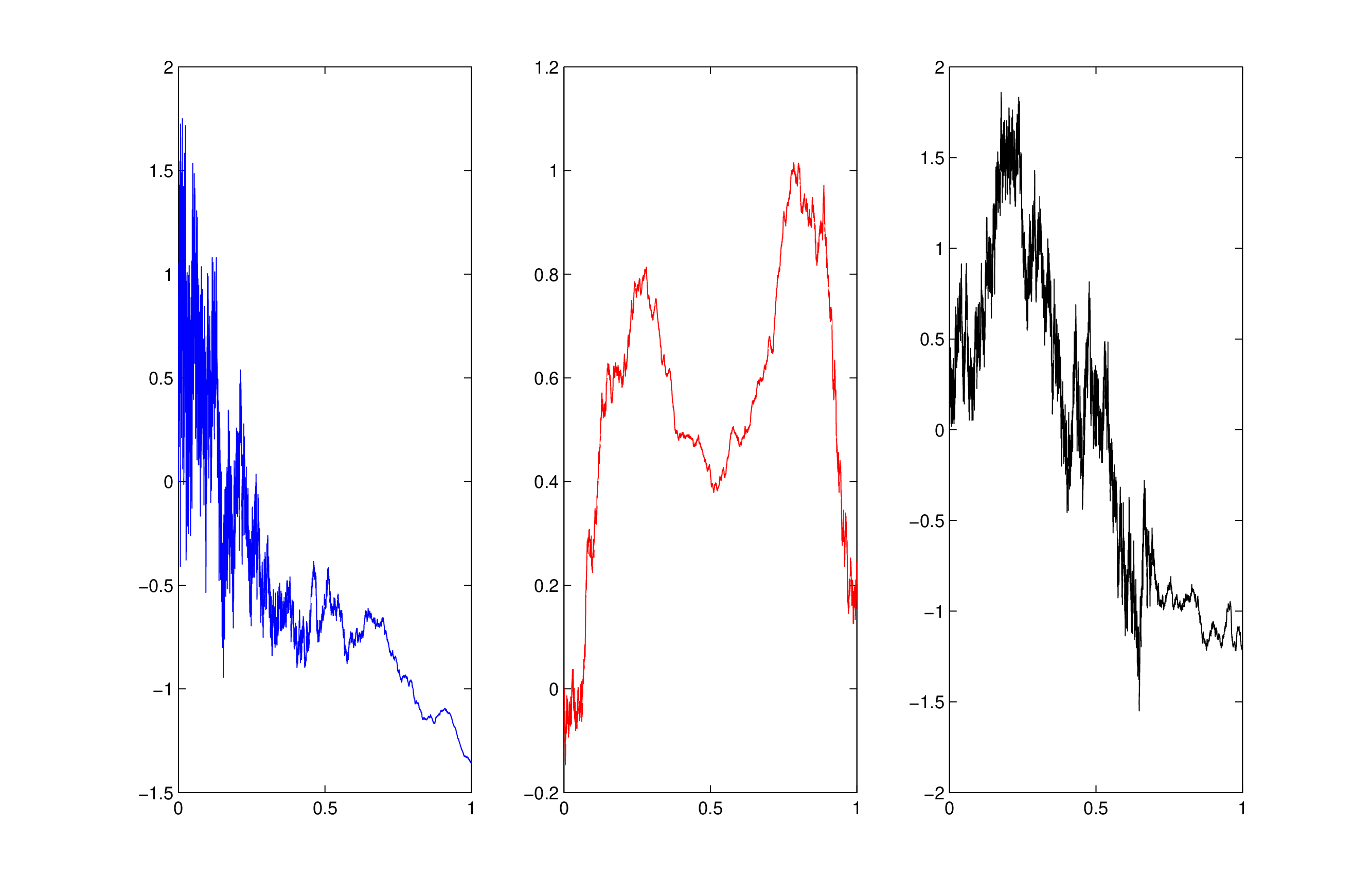}
\end{tabular}
\caption{\emph{Simulated mBm process with $H(\cdot)$ linear
function (left), $H(\cdot)$ periodic function (middle) and
$H(\cdot)$ logistic function (right).}} \label{fig6}
\end{center}
\end{figure}
Then, for each sample $\mathbf{(\mathcal{C}_i)}$ with $i \in
[4,6]$, we have estimated the Hurst function $\widehat{H}_{n,i}$
by using the localized version of IRS with $\gamma=0.3$ and the
GQV method. We note that both methods
provide correct results represented by Figure~\ref{fig7} below.\\

\begin{figure}[htbp]
\begin{center}
\begin{tabular}{c}
\includegraphics[width=15cm,height=5cm]{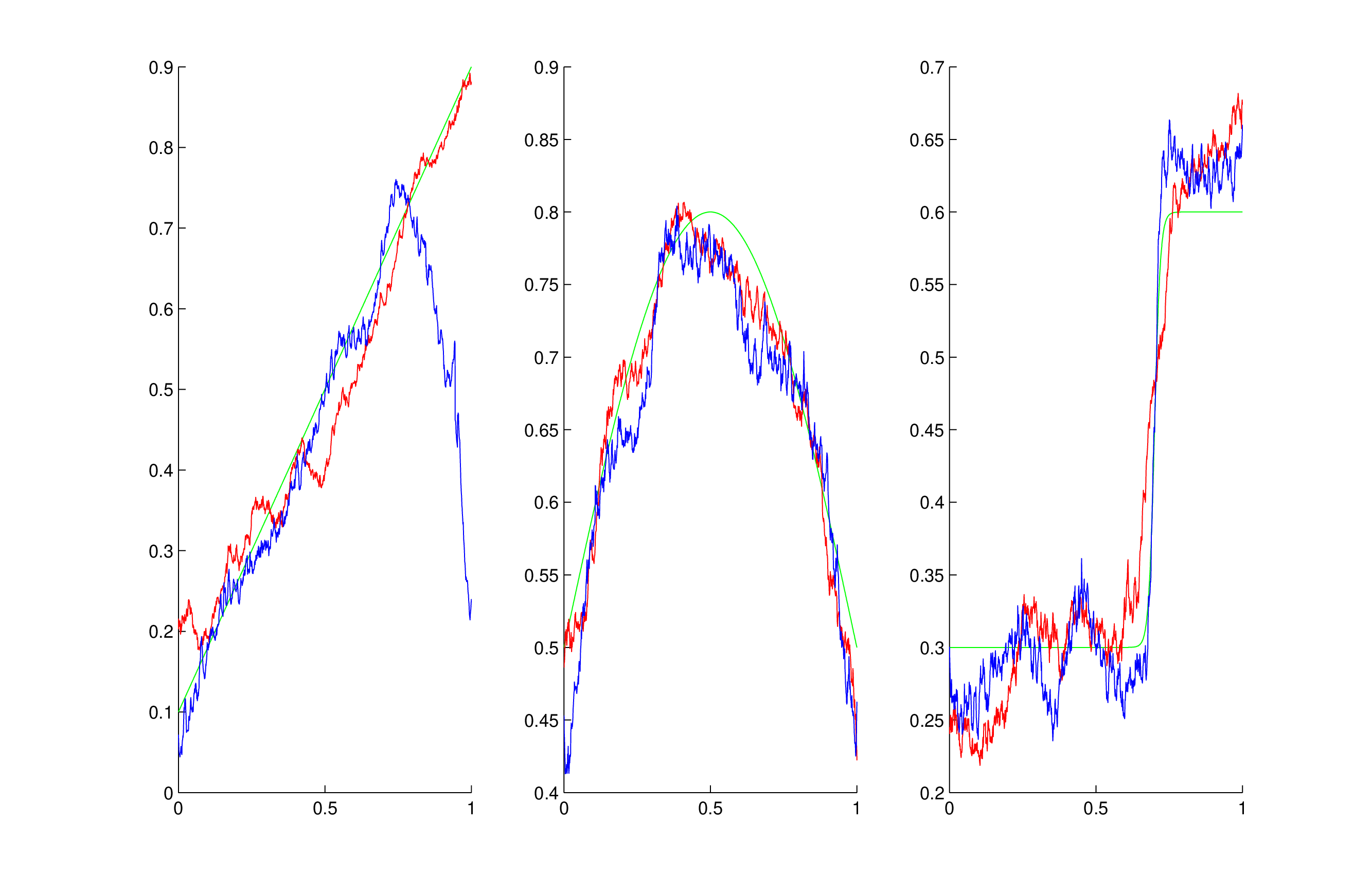}
\end{tabular}
\caption{\emph{Estimation of the Hurst function $H(\cdot)$ with
$H(\cdot)$ linear function (left), $H(\cdot)$ periodic function
(middle) and $H(\cdot)$ logistic function (right). The graphs of
function $H(\cdot)$, $\widehat{H}_n^{(IRS)}(\cdot)$ its estimation
by IRS, and $\widehat{H}_n^{(GQV)}(\cdot)$ its estimation by GQV,
are represented in green, red and blue respectively.}}
\label{fig7}
\end{center}
\end{figure}
These results are plainly confirmed by Monte Carlo simulations.
Actually, for each case $\mathbf{(\mathcal{C}_i)}$ with $i \in
[4,6]$, we have made $M=1000$ simulations of independent copies of
mBm sequences $B^{(k)}_{H_i(\cdot)}=\left(B^{(k)}_{H_i(t_0)}(t_0),
\ldots, B^{(k)}_{H_i(t_n)}(t_n)\right)$, for $k=1,\ldots,M$. Then
we have computed the Mean Integrate Square Error (MISE) defined as
$\ds \text{MISE}=\E \left(\frac{1}{n+1} \sum_{j=0}^n
\left|\widehat{H}(t_j) - H(t_j)\right|^2 \right) = \E \|
\widehat{H} - H \|_{L^2(0,1)}^{2}$ which is a criterion widely
used in functional estimation, see Table~\ref{tab3} below.\\
\text{}
\begin{table}[htbp]
\begin{center}
\begin{tabular}{|c|c|c|c|}
  \hline
   \  & $H(\cdot)$ Linear   &  $H(\cdot)$ Periodic & $H(\cdot)$ Logistic\\
  \hline
  MISE by IRS & $2.6743 \times 10^{-4}$ & $1.4743 \times 10^{-4}$ &$5.3546 \times 10^{-3}$ \\
  \hline
  MISE by GQV & $8.9547 \times 10^{-4}$  & $5.4743 \times 10^{-4}$& $8.9743 \times 10^{-4}$ \\
  \hline
\end{tabular}
\caption{MISE given by IRS method and GQV method} \label{tab3}
\end{center}
\end{table}
We observe through Table~\ref{tab3} that both methods provide
globally the same results when the function $H(\cdot)$ varies
slowly (see linear and periodic cases), whereas in the case where
$H(\cdot)$ presents the abrupt variation it appears that the GQV
is a bit more precise compared to the IRS method.

\section{Proofs of the main results}

This section contains the proof of the results of Section 3. Note
that we have divided the proof of Theorem 1 in two parts: first we
consider the general case of Gaussian processes with stationary
increments and then in a second part we investigate the
application to fractional Bronwnian motion.

\subsection{Proof of localization}

First, we can deduce as a corollary that the
$a$-Generalized increments sequence $\left ( \Delta_a X(t_k)
\right )_{0 \leq k \leq n-L-1}$ form a family of stationary
identically distributed centered Gaussian r.v. with variance

\ban \sigma^2_{a,n} &=& \cov \left( \Delta_a X(t_{k}) , \Delta_a
X(t_{k}) \right)\\&=& \int_{\R} |g_a(\xi/n)|^2 \cdot f(\xi) \;
d\xi \\&=& 2 \int_{\R_{+}} |g_a(\xi/n)|^2 \cdot f(\xi) \; d\xi,,
\ean covariance given, for all $0 \leq k_1, k_2 \leq n-L-1$, by
\ban \cov_{a,n}(t_{k_1},t_{k_2}) &=& \cov \left( \Delta_a
X(t_{k_1}) , \Delta_a X(t_{k_2}) \right)\\&=&\int_{\R} e^{i
(k_1-k_2)\xi/n} |g_a(\xi/n)|^2 \cdot f(\xi) \; d\xi\\ &=& 2
\int_{\R_{+}} \cos \left((k_1-k_2)\xi/n\right) |g_a(\xi/n)|^2
\cdot f(\xi) \; d\xi, \ean and correlation between two successive
$a$-Generalized increments defined by \ban \rho_{a,n} = \frac{\cov
\left( \Delta_a X(t_{k+1}) , \Delta_a X(t_{k}) \right)}{\left [
\cov \left( \Delta_a X(t_{k+1}) , \Delta_p X(t_{k+1})
\right)\right]^{1/2} \cdot \left [ \cov \left( \Delta_a X(t_{k}) ,
\Delta_p X(t_{k})
\right)\right]^{1/2}}=\frac{\cov_{a,n}(t_{k+1},t_k)}{\sigma^2_{a,n}},
\label{rho_an} \ean where $g_a(\cdot)$ is described by
$\eqref{g_a}$. Therefore, for a fixed $0 \leq k \leq n-L-1$, it is
easy to remark that there exist two independent standard Gaussian
r.v. $Z_{k}, Z_{k+1} \egalelaw \mathcal{N}(0,1)$ such that
 \ba
 \Delta_a X(t_k)&=& \sigma_{a,n} Z_{k} \label{d1} \\
 \Delta_a X(t_{k+1})&=& \sigma_{a,n} \left( \rho_{a,n}  Z_{k} +
 \sqrt{1-\rho^2_{a,n}} Z_{k+1}\right), \label{d2} \ea
where the sign $\egalelaw$ means equal in distribution.

\remark In the particular case of the fBm, the correlation between
two successive $a$-Generalized increments, $\rho_{a,n}$, does not
depend on $n$. Indeed, we know that the spectral density of the
fBm is given by $\eqref{densiteSpec}$, then we have \ban
\rho_{a,n}(H)= \frac{\int_{\R_{+}} \cos(\xi/n) |g_a(\xi/n)|^2
\cdot \xi^{-(2H+1)} \; d\xi}{\int_{\R_{+}} |g_a(\xi/n)|^2 \cdot
\xi^{-(2H+1)} \; d\xi}. \ean And after, we can change variable
$\xi/n$ to $u$. So this implies that \ban \rho_{a,n}=
\frac{\int_{\R_{+}} \cos(u) |g_a(u)|^2 \cdot u^{-(2H+1)} \; du}{
\int_{\R_{+}} |g_a(u)|^2 \cdot u^{-(2H+1)} \; du}=\rho_a(H), \ean
which is independent of $n$.

\subsection{Proof of CLT for Gaussian Processes with stationary increments}
The proof uses the notion
of Hermite rank and Breuer-Major theorem, see for e.g Arcones
\cite[Theorem 4, p.2256]{Arcones:1994} or Nourdin \emph{et al}
\cite[Theorem 1, p.2]{Nourdin:etal:2010}.

\begin{defi}[Hermite rank] \label{def:definition} Let $G$ be a $\R^d$ Gaussian vector and
$\phi:\R^d \rightarrow \R$ be a measurable function such that $\ds
\E \left| \phi(G) \right|^2 < +\infty$. Then, the function $\phi$
is said to have Hermite rank equal to the integer $q \geq 1$ with
respect to Gaussian vector $G$, if (a) $\ds \E \left[
\left(\phi(G) - \E \left( \phi(G) \right) \right)P_m(G) \right]=0$
for every polynomial $P_m$ $\left(\text{on } \R^d\right)$ of
degree $m \leq q-1$; and (b) there exists a polynomial $P_q$
$\left(\text{on } \R^d\right)$ of degree $q$ such that $\ds \E
\left[ \left(\phi(G) - \E \left( \phi(G) \right) \right)P_q(G)
\right] \neq 0$.
\end{defi}

We first give the proof of Theorem 1 in the general case of
Gaussian processes with stationary increments and then in a
separate part the application to fractional Brownian motion.

\begin{prfT}
First, in the sequel we denote by $G=\left( G_{k},G_{k+1} \right)=
\left( \Delta_a X(t_k), \Delta_a X(t_{k+1})\right )$ these two
successive stationary $a$-Generalized increments defined by
$(\eqref{d1},\eqref{d2})$. Then, according to Bardet and Surgailis
\cite[Appendix, p.31]{Bardet:Surgailis:2009}, we know that
 \ban
 \E  \left [ \psi\left( G_{k}, G_{k+1} \right ) \right] = \Lambda_0(\rho_{a}), &\text{ and }&
 \E  \left | \psi\left( G_{k}, G_{k+1} \right ) \right
 |^2 < +\infty,
  \ean
where $\Lambda_0(.)$ is defined by $\eqref{L0}$ and $\rho_{a}$ is
the correlation between $G_{k}$ and $G_{k+1}$. To achieve our
goal, we start by defining a new function $\phi:R^2 \rightarrow
\R$ such that \ban \phi(X,Y)=\psi(X,Y)-\Lambda_0(\rho_{a}). \ean
Then, $\phi$ is in fact a Hermite function with respect to
Gaussian vector $G=\left( G_{k},G_{k+1} \right)$ with rank equal
to 2. Therefore, by applying Breuer-Major theorem, see e.g Arcones
\cite[Theorem 4, p.2256]{Arcones:1994} or Nourdin \emph{et al}
\cite[Theorem 1, p.2]{Nourdin:etal:2010}, we get directly the CLT
$\eqref{TCLGauss}$. So, the key argument of our proof is to
determine the Hermite rank of $\phi$. We include here the proof of
the fact that the Hermite rank is 2 as the proof does not seem to
appear elsewhere. Let $P_0(X,Y)=c_0,
P_1(X,Y)=c_{11}X+c_{12}Y+c_{10}$ and $P_2(X,Y)=X^2$ be three
polynomials (on $\R^2$) with degree respectively 0, 1 and 2.
First, it is easy to see that $\ds \E \left[ \phi(G) P_0(G)
\right]=0$. Now, we must to show that $\ds \E \left[ \phi(G)
P_1(G) \right]=0$. We have \ban  \E \left[ \phi(G) P_1(G) \right]=
c_{11} \E \left[ \phi(G) G_{k} \right] + c_{12} \E \left[ \phi(G)
G_{k+1} \right] +c_{10} \underbrace{\E \left[ \phi(G)
\right]}_{=0}. \ean Then, \ban  \E \left[ \phi(G) P_1(G) \right]=
c_{11}\E \left[ \psi(G_{k},G_{k+1}) G_{k} \right]-
c_{11}\Lambda_0(\rho_{a})\underbrace{\E[G_k]}_{=0}+c_{12}\E \left[
\psi(G_{k},G_{k+1}) G_{k+1} \right]-
c_{12}\Lambda_0(\rho_{a})\underbrace{\E[G_{k+1}]}_{=0}, \ean
because $G_{k}$ and $G_{k+1}$ are zero-mean r.v. and due to the
fact that the r.v $G_k$ and $G_{k+1}$ have a symmetric function,
we can write without any restrictions that \ban  \E \left[ \phi(G)
P_1(G) \right]= (c_{11} + c_{12}) \E \left[ \psi(G_{k},G_{k+1})
G_{k} \right]. \ean By using definition of $\left( G_k,G_{k+1}
\right)= \left( \Delta_a X(t_k), \Delta_a X(t_{k+1})\right )$
given by $(\eqref{d1},\eqref{d2})$, we get \ban (c_{11} +
c_{12})^{-1} \E \left[ \phi(G) P_1(G) \right]= \sigma_{a} \E
\left[ \psi \left(\sigma_{a} Z_k,\sigma_{a} \left( \rho_{a} Z_k +
\sqrt{1-\rho_{a}^2} Z_{k+1}\right)\right) Z_{k} \right] \ean where
$Z_k$ and $Z_{k+1}$ are two independent standard Gaussian r.v.
$Z_{k}, Z_{k+1} \egalelaw \mathcal{N}(0,1)$. Thus, by using
homogeneity property of $\psi(\cdot,\cdot)$ specified by:
$\psi(aX,aY)=\psi(X,Y)$, we obtain \ban (c_{11} + c_{12})^{-1} \E
\left[ \phi(G) P_1(G) \right]= \sigma_{a} \E \left[ \psi
\left(Z_k,\left( \rho_{a} Z_k + \sqrt{1-\rho_{a}^2}
Z_{k+1}\right)\right) Z_{k} \right]. \ean Next, we have \ban
(c_{11} + c_{12})^{-1} \E \left[ \phi(G) P_1(G) \right] =
\frac{\sigma_{a,n}}{2\pi} \int_{\R^2} \psi \left(z_1, \rho_{a} z_1
+ \sqrt{1-\rho_{a}^2} z_2 \right) z_1
\exp\left(-\frac{z_1^2+z_2^2}{2}\right) dz_1 dz_2. \ean And after,
we can change variables $(z_1,z_2)$ to polar coordinates $(r \cos
(\theta),r \sin(\theta))$. So, this implies that \ban  2\pi
\sigma_p^{-1}(c_{11} + c_{12})^{-1} \E \left[ \phi(G) P_1(G)
\right] &=&  \int_{-\pi}^{\pi} \psi \left(\cos(\theta), \rho_{a}
\cos(\theta) + \sqrt{1-\rho_{a}^2} \sin(\theta) \right) \cos(\theta) d \theta  \\
&\times& \underbrace{\int_0^{+\infty} r^2 \exp\left(-\frac{r^2}{2}
\right)dr}_{=\sqrt{\pi/2}}. \ean  Then, we remark that \ban
\int_{-\pi}^{0}  \psi \left(\cos(\theta), \rho_{a} \cos(\theta)
+
\sqrt{1-\rho_{a}^2} \sin(\theta) \right) \cos(\theta) d\theta &=& \\
\int_{-\pi}^{0} \frac{ \left |\cos(\theta) +  \rho_{a}
\cos(\theta) + \sqrt{1-\rho_{a}^2} \sin(\theta)\right| }{\left
|\cos(\theta)\right| +  \left |\rho_{a} \cos(\theta) +
\sqrt{1-\rho_{a}^2} \sin(\theta)\right|} \cos(\theta) d\theta &\underbrace{=}_{u=\theta+\pi}&\\
\int_{0}^{\pi} \frac{ \left |\cos(u-\pi) +  \rho_{a} \cos(u-\pi)
+ \sqrt{1-\rho_{a}^2} \sin(u-\pi)\right| }{\left
|\cos(u-\pi)\right| + \left |\rho_{a} \cos(u-\pi) +
\sqrt{1-\rho_{a}^2} \sin(u-\pi)\right|} \cos(u-\pi) du &=& \\ -
\int_{0}^{\pi} \frac{ \left |\cos(u) + \rho_{a} \cos(u) +
\sqrt{1-\rho_{a}^2} \sin(u)\right| }{\left |\cos(u)\right| +
\left |\rho_{a} \cos(u) + \sqrt{1-\rho_{a}^2} \sin(u)\right|}
\cos(u) du &=&
\\ - \int_{0}^{\pi}  \psi \left(\cos(\theta), \rho_{a} \cos(\theta)
+ \sqrt{1-\rho_{a,n}^2} \sin(\theta) \right) \cos(\theta) d\theta.
\ean Thus, we deduce directly that $\ds \E \left[ \phi(G) P_1(G)
\right]=0$. In the similar way, it is easy to prove that $\ds \E
\left[ \phi(G) P_2(G) \right]\neq0$. So, by using
Definition~\ref{def:definition}, we can say that $\phi$ is a
Hermite function with rank equal to 2. Therefore,
Theorem~\ref{th:Theorem1} becomes an application of Breuer-Major
theorem which gives directly the proof of CLT $\eqref{TCLGauss}$.
\hfill$\blacksquare$ \end{prfT}

\subsection*{Proof of CLT for FBm}
Next, we present the correlation function properties of the
$a$-Generalized increments sequence of a fBm.
\begin{prop}[Correlation function of the
$a$-Generalized increments] \label{prop:Property} Let $\left
(B_H(t), t \in [0,1] \right)$ be a fBm with Hurst parameter $H \in
\left (0,1 \right)$ and let $\left ( \Delta_a B_H(t_k) \right )_{0
\leq k \leq n-L-1}$ its $a$-Generalized increments sequence
defined by $\eqref{delta:tk}$, with $a \in \mathcal{A}(p,L)$ a
filter given by $\eqref{filter}$. Then, for all $j \in \ZZ$, we
have \ban r_{a,n}(j) = -\frac{1}{2 \cdot n^{2H}} \times C_a(j)
\ean where $r_{a,n}(j) := \cov \left( \Delta_a B_H(t_0), \Delta_a
B_H(t_j) \right)$ and $C_a(j)$ is given by \ban C_a(j) &:=&
\sum_{l_1,l_2=0}^L a_{l_1}a_{l_2} |j+l_2-l_1|^{2H}
\\ &\underset{j\rightarrow +\infty }{\sim
}&  \binom{2H}{2p} \cdot\left (\sum_{l=0}^{L} a_l l^p \right)^2
\times j^{2H-2p} \text{  with  } \binom{2H}{2p}= \ds
\frac{\prod_{k=0}^{2p-1}(2H-k)}{(2p)!}.  \ean And the correlation
between two successive $a$-Generalized increments, is specified by
\ba \label{rho_aH}  \rho_{a}(H) = \frac{ \sum_{l_1,l_2=0}^L
a_{l_1}a_{l_2}|1+l_2-l_1|^{2H}}{\sum_{l_1,l_2=0}^L
a_{l_1}a_{l_2}|l_2-l_1|^{2H}}. \ea
\hfill $\square$
\end{prop}

\begin{prfP}
To compute the covariance function of the $a$-Generalized
increments sequence, we start by using the initial formula of the
covariance function of a fBm defined by $\eqref{covFBM}$. Then, we
obtain \ban r_{a,n}(j) &:=& \cov \left( \Delta_a B_H(t_0),
\Delta_a B_H(t_j) \right) \\ &=& \sum_{l_1,l_2=0}^L a_{l_1}
a_{l_2} \E \left [B_H(t_{l_1})B_H(t_{j+l_2}) \right] \\ &=&
\frac{1}{2} \underbrace{\sum_{l_1,l_2=0}^L a_{l_1} a_{l_2}
t_{l_1}^{2H}}_{=0}+ \frac{1}{2} \underbrace{\sum_{l_1,l_2=0}^L
a_{l_1} a_{l_2} t_{j+l_2}^{2H}}_{=0} - \frac{1}{2}
\sum_{l_1,l_2=0}^L a_{l_1} a_{l_2} |t_{j+l_2}-t_{l_1}|^{2H} \\ &=&
-\frac{1}{2 \cdot n^{2H}} \times C_a(j) \ean where $ C_a(j) =
\sum_{l_1,l_2=0}^L a_{l_1}a_{l_2} |j+l_2-l_1|^{2H} $. Now, we give
an equivalent of $C_a(j)$ when $j \rightarrow + \infty$. To do
this, we use the Taylor expansion as follows \ban C_a(j)&=& j^{2H}
\times \sum_{l_1,l_2=0}^L a_{l_1} a_{l_2} \left |
1+\frac{l_1-l_2}{j} \right|^{2H} \\ &=& j^{2H} \times \left [
\underbrace{\sum_{l_1,l_2=0}^L a_{l_1} a_{l_2} \cdot 1}_{=0} +
\sum_{k=1}^{+\infty} \binom{2H}{k}j^{-k} \sum_{l_1,l_2=0}^L
a_{l_1} a_{l_2}|l_2-l_1|^k \right ]. \ean Next, by using
$\eqref{doubleSum}$, we know that when we sum over $k$, every term
in the expansion gives a zero contribution for any integer $k <
2p$. So this implies that \ban C_a(j) &=&  \binom{2H}{2p}
\cdot \left (\sum_{l=0}^{L} a_l l^p \right)^2 \times j^{2H-2p} + o \left( j^{2H-2p} \right) \\
&\underset{j\rightarrow +\infty }{\sim }& \binom{2H}{2p} \cdot
\left (\sum_{l=0}^{L} a_l l^p \right)^2 \times j^{2H-2p}. \ean
This finishes the proof of Property~\ref{prop:Property}.
\hfill$\blacksquare$
\end{prfP}
\noindent And after,
we note that the function
$\psi(\cdot,\cdot)$ satisfies the homogeneity property specified
by: $\psi(aX,aY)=\psi(X,Y)$. So, this allows us to rewrite
$\IRS_{a,n}(B_H)$ as follows \ban \IRS_{a,n}(B_H)= \frac{1}{n-L}
\sum_{k=0}^{n-L-1} \psi \left( \Delta_a B_H^{\text{std}}(t_k) ,
\Delta_a B_H^{\text{std}}(t_{k+1}) \right) \ean where $\Delta_a
B_H^{\text{std}}$ represents the standardized version of $\Delta_a
B_H$ described, for all $0 \leq k \leq n-L-1$, as \ban \Delta_a
B_H^{\text{std}}(t_k)=\frac{\Delta_a B_H(t_k)}{\sqrt{Var\left[
\Delta_a B_H(t_k) \right]}} \ean and its covariance function is
given by \ban r_{a}(j)&:=&\cov \left( \Delta_a
B_H^{\text{std}}(t_0), \Delta_a B_H^{\text{std}}(t_j) \right) \\
&=& \frac{r_{a,n}(j)}{r_{a,n}(j)}
\\ &=& \frac{C_{a}(j)}{C_{a}(0)}  \\
&\underset{j\rightarrow +\infty }{\sim }& \ds \frac{1}{C_a(0)}
\cdot \binom{2H}{2p} \cdot \left (\sum_{l=0}^{L} a_l l^p \right)^2
\times j^{2H-2p} \ean which is independent of $n$. So, according
to Theorem~\ref{th:Theorem1},  the key argument is to prove that
\ban \sum_{j \in \ZZ} |r_a(j)|^2 <+\infty.\ean Thus, by using a
Riemman sum argument, we can deduce immediately that this
condition is verified if and only if $4H - 4p < -1$, i.e
$H<p-1/4$, and this implies that $H \in (0,3/4)$ if $p=1$ and that
$H \in (0,1)$ if $p=2$. Therefore, the assumption $\eqref{condG}$
of Theorem~\ref{th:Theorem1} is satisfied and so we obtain a
simple intuitive proof of the CLT $\eqref{TCL:fBm}$ applied to the
IRS of fBm. This finishes the proof of
Theorem~\ref{th:Theorem1}.\hfill$\blacksquare$ \text{}\\
\noindent {\bf End of The proof of Theorem 1}\\
\subsection{Proof of CLT for mBm}
The proof of Theorem 2 relies on a localization argument given in the following Lemma

\begin{lemm}
\label{lemm:Lemma} First, we consider $t^{*} \in (0,1)$ be an
arbitrary fixed point, $\gamma \ (0,1)$ a fixed parameter which
allows to control the size of the indices set around $t^{*}$, and
$a \in \mathcal{B}(p)$
a binomial filter. \\
Let $B_{H(.)}=(B_{H(t)}(t), t \in [0,1])$ be a mBm with Hurst
function $H(\cdot) \in \mathcal{C}^{\eta} \left ([0,1],
\left[H_{\diamond}, H^{\diamond}\right ] \right)$ and
$B_{H^{*}}=(B_{H^{*}}(t), t \in [0,1])$ a fBm with Hurst index
$H^{*}=H(t^*)$. \\
Moreover, we consider
$\IRS_{a,n}^{\gamma,t^{*}}\left(B_{H(\cdot)}\right)$ the localized
version of IRS for mBm defined by $\eqref{IRSmBm}$, and
$\IRS_{a,n}\left(B_{H^{*}}\right)$ a modified version of the IRS
for fBm described as follows \ban
\IRS_{a,n}\left(B_{H^{*}}\right)=\frac{1}{ 2n^{1-\gamma} +1}
\sum_{k=\lfloor nt^{*}-n^{1-\gamma} \rfloor}^{\lfloor
nt^{*}+n^{1-\gamma} \rfloor} \psi \left ( \Delta_a B_{H^{*}}(t_k),
\Delta_a B_{H^{*}}(t_{k+1}) \right ). \ean  Then \ba
\label{delta:IRS} \E \left |
\IRS_{a,n}^{\gamma,t^{*}}\left(B_{H(\cdot)}\right) -
\IRS_{a,n}\left(B_{H^{*}}\right) \right |^2 =
\underset{n\rightarrow +\infty }{\mathcal{O}} \left( n^{-\gamma
\eta} \right ).\ea
\end{lemm}

\begin{prfL}
For n large enough, we have \ban \left | \Delta
\IRS_{n}^{\gamma,\eta} \right |^2 &:=& \left |
\IRS_{a,n}^{\gamma,t^{*}}\left(B_{H(\cdot)}\right) -
\IRS_{a,n}\left(B_{H^{*}}\right) \right |^2 \\ &=& \frac{1}{
2n^{1-\gamma} +1}  \left | \sum_{k=\lfloor nt^{*}-n^{1-\gamma}
\rfloor}^{\lfloor nt^{*}+n^{1-\gamma} \rfloor}  \psi \left (
\Delta_a B_{H(t_k)}(t_k), \Delta_a B_{H(t_{k+1})}(t_{k+1}) \right)
- \psi \left ( \Delta_a B_{H^{*}}(t_k), \Delta_a
B_{H^{*}}(t_{k+1}) \right ) \right |^2. \ean Then, by using
Cauchy-Schwarz inequality, we get \ban \left | \Delta
\IRS_{n}^{\gamma,\eta} \right |^2 \leq \frac{1}{ 2n^{1-\gamma} +1}
\sum_{k=\lfloor nt^{*}-n^{1-\gamma} \rfloor}^{\lfloor
nt^{*}+n^{1-\gamma} \rfloor} \left | \psi \left (\Delta_a
B_{H(t_k)}(t_k),\Delta_a B_{H(t_{k+1})}(t_{k+1}) \right) - \psi
\left (\Delta_a B_{H^{*}}(t_k), \Delta_a B_{H^{*}}(t_{k+1}) \right
) \right |^2. \ean This implies that, \ban \E \left | \Delta
\IRS_{n}^{\gamma,\eta} \right |^2 \leq \frac{1}{ 2n^{1-\gamma} +1}
\sum_{k=\lfloor nt^{*}-n^{1-\gamma} \rfloor}^{\lfloor
nt^{*}+n^{1-\gamma} \rfloor} \E \left | \psi \left (\Delta_a
B_{H(t_k)}(t_k), \Delta_a B_{H(t_{k+1})}(t_{k+1}) \right) - \psi
\left (\Delta_a B_{H^{*}}(t_k),\Delta_a B_{H^{*}}(t_{k+1}) \right
) \right |^2. \ean Now, we recall that $\Omega^{*}$ represents the
event with probability 1 introduced in Subsection 2.3. Then,
according to Bru\v{z}ait\.{e} \& Vai\v{c}iulis \cite[Lemma~1,
formula 3.3, p. 262]{Bruzaite:Vaiciulis:2008} and our Taylor
expansion $\eqref{taylor}$, we deduce that \ban \E \left \{ \left
| \psi \left ( \Delta_a B_{H(t_k)}(t_k),\Delta_a
B_{H(t_{k+1})}(t_{k+1}) \right) - \psi \left (\Delta_a
B_{H^{*}}(t_k), \Delta_a
B_{H^{*}}(t_{k+1}) \right ) \right |^2 \I_{\Omega^{*}} \right\}\\
\leq C\left( \rho_a(H^*) \right) \cdot \E \left \{ \left (
\left|\Delta_a R(t_{k}) \right|+ \left|\Delta_a R(t_{k+1}) \right|
\right ) \I_{\Omega^{*}} \right\} .\ean where the constant
$C\left( \rho_a(H^*) \right)>0$ depend only on $\rho_a(H^*)$ and
$\left (\Delta_a R(t_{k})\right)$ corresponds to $a$-Generalized
increments at $t_k$ of the rest $(R(t), t \in [0,1])$ defined by
$\eqref{rest}$. Next, by using $\eqref{supR}$, we deduce that
there exist a constant $\kappa > 0$ such as \ban \E \left | \psi
\left ( \Delta_a B_{H(t_k)}(t_k),\Delta_a B_{H(t_{k+1})}(t_{k+1})
\right) - \psi \left (\Delta_a B_{H^{*}}(t_k), \Delta_a
B_{H^{*}}(t_{k+1}) \right ) \right |^2 \I_{\Omega^{*}} \leq \kappa
n^{-\gamma \eta} . \ean Therefore, we obtain \ban \E \left \{
\left | \Delta \IRS_{n}^{\gamma,\eta} \right |^2
\I_{\Omega^{*}}\right\} \leq \kappa n^{-\gamma \eta} . \ean
Moreover, we know that \ban \E \left | \Delta
\IRS_{n}^{\gamma,\eta} \right |^2 = \E \left \{ \left | \Delta
\IRS_{n}^{\gamma,\eta} \right |^2
\I_{\Omega^{*}}\right\}+\underbrace{\E \left \{ \left | \Delta
\IRS_{n}^{\gamma,\eta} \right |^2 \I_{\Omega \backslash
\Omega^{*}}\right\}}_{=0},\ean by applying Cauchy-Schwarz
inequality. Then, we get directly \ban \E \left | \Delta
\IRS_{n}^{\gamma,\eta} \right |^2 \leq \kappa n^{-\gamma \eta}
.\ean This finishes the proof of Lemma~\ref{lemm:Lemma}.
\hfill$\blacksquare$

 \begin{prfT}
First, according to our Lemma~\ref{lemm:Lemma}, we have \ban \E
\left | \IRS_{a,n}^{\gamma,t^{*}}\left(B_{H(\cdot)}\right) -
\IRS_{a,n}\left(B_{H^{*}}\right) \right |^2 =
\underset{n\rightarrow +\infty }{\mathcal{O}} \left( n^{-\gamma
\eta} \right ). \ean Next, it is easy to see that \ban \E \left |
n^{(1-\gamma)/2} \left(
\IRS_{a,n}^{\gamma,t^{*}}\left(B_{H(\cdot)}\right)-\Lambda_a(H^*)
\right) - n^{(1-\gamma)/2} \left(
\IRS_{a,n}\left(B_{H^{*}}\right)-\Lambda_a(H^*) \right) \right |^2
= \underset{n\rightarrow +\infty }{\mathcal{O}} \left(
n^{1-\gamma(1+ \eta)} \right ). \ean After, by applying
Theorem~\ref{th:Theorem1}, we know that \ban n^{(1-\gamma)/2}
\left( \IRS_{a,n}\left(B_{H^{*}}\right)-\Lambda_a(H^*) \right)
\convlaw \mathcal{N}(0,\Sigma^2_a) \text{ with } \left\{
\begin{array}{ll}
    \ H^* \in \left (0,3/4 \right) & \mbox{if } p=1 \\
    \ H^* \in \left (0, 1 \right) & \mbox{if } p \geq 2 \\
\end{array}
\right.\ean Therefore, CLT $\eqref{TCLmBm}$ is satisfied if and
only if $\gamma(1+\eta)>1$.\\
Let us now sketch how to extend it to the multidimensional case: first we may operate a multidimensional freezing of time in the sense that there exists an almost sure event $\Omega^*$ such that
$$\forall i,\, B(H(t),t)1_{\Omega^*}=B^i(H(t_i,t)1_{\Omega^*}+R^i(t)1_{\Omega^*},\qquad  \sup_{s
\in [0,1]} \left |R^i(s) \I_{\Omega^{*}} \right | \leq C^{*}(\omega)
|H(t)-H^{*}| $$
and the process $B^i$ are defined using wavelet expansion so that the correlations between them are well described. We may then consider fractional Brownian motions rather than mBm. Secondly we use Cramer-Wold device (see e.g. Th. 7.7 in Billingsley \cite{Billingsley:1986}) : it is sufficient to get the CLT for every real numbers $b_1,...,b_m$ for
$$ n^{(1-\gamma)/2}\sum_{i=1}^m b_i\,\left(\IRS^{\gamma,t_i}_{a,n}(B^i_{H(t_i)})-\Lambda_a(H(t_i))\right)$$
which is obtained exactly as before.
\hfill$\blacksquare$
\end{prfT}

\end{prfL}


\bibliographystyle{acm}
\bibliography{IRSmBm}

\begin{thebibliography}{10}

\bibitem{Abry:etal:2003}
{\sc Abry, P., Flandrin, P., Taqqu, M.~S., and Veitch, D.}
\newblock {\em Self-similarity and long-range dependence through the wavelet
  lens, in Theory and applications of long-range dependenc}.
\newblock Birkhauser, Boston., 2003.

\bibitem{Arcones:1994}
{\sc Arcones, M.~A.}
\newblock Limit theorems for nonlinear functionals of a stationary gaussian
  sequence of vectors.
\newblock {\em Ann. Probab 22\/} (1994), 2242--2274.

\bibitem{ayache:etal:2007}
{\sc Ayache, A., Bertrand, P., and L{\'e}vy~V{\'e}hel, J.}
\newblock A central limit theorem for the generalized quadratic variation of
  the step fractional {B}rownian motion.
\newblock {\em Stat. Inference Stoch. Process. 10}, 1 (2007), 1--27.

\bibitem{Bardet:Bertrand:2007}
{\sc Bardet, J.-M., and Bertrand, P.}
\newblock Identification of the multiscale fractional {B}rownian motion with
  biomechanical applications.
\newblock {\em J. Time Ser. Anal. 28}, 1 (2007), 1--52.

\bibitem{Bardet:Bertrand:2007b}
{\sc Bardet, J.~M., and Bertrand, P.~R.}
\newblock Definition, properties and wavelet analysis of multiscale fractional
  brownian motion.
\newblock {\em Fractals 15\/} (2007), 73--87.

\bibitem{Bardet:Bertrand:2010}
{\sc Bardet, J.~M., and Bertrand, P.~R.}
\newblock A nonparametric estimator of the spectral density of a
  continuous-time gaussian process observed at random times.
\newblock {\em Scandinavian Journal of Statistics\/} (2010), 1--41.

\bibitem{Bardet:Surgailis:2009}
{\sc Bardet, J.~M., and Surgailis, D.}
\newblock Measuring roughness of random paths by increment ratios.
\newblock {\em HAL : hal-00238556, version 2\/} (2009).

\bibitem{Benassi:etal:1998}
{\sc Benassi, A., Cohen, S., and Istas, J.}
\newblock Identifying the multifractional function of a gaussian process.
\newblock {\em Statist. Probab. Lett. 39\/} (1998), 337--345.

\bibitem{Benassi:etal:1997}
{\sc Benassi, A., Jaffard, S., and Roux, D.}
\newblock Gaussian processes and pseudodifferential elliptic operators.
\newblock {\em Rev. Mat. Iberoam. 13\/} (1997), 19--81.

\bibitem{Bertrand:etal:2010}
{\sc Bertrand, P.~R., Hamdouni, A., and Khadhraoui, S.}
\newblock Modelling nasdaq series by sparse multifractional brownian motion.
\newblock {\em Methodology and Computing in Applied Probability\/} (2010).

\bibitem{Billingsley:1986}
{\sc Billingsley, P.}
\newblock {\em Probability and measure}, second~ed.
\newblock Wiley Series in Probability and Mathematical Statistics: Probability
  and Mathematical Statistics. John Wiley \& Sons Inc., New York, 1986.

\bibitem{Bruzaite:Vaiciulis:2008}
{\sc Bru\v{z}ait\.{e}, K., and Vai\v{c}iulis, M.}
\newblock The increment ratio statistic under deterministic trends.
\newblock {\em Lithuanian Mathematical Journal 48\/} (2008), 256--269.

\bibitem{Chen:Wood:1998}
{\sc Chen, G., and Wood, A. T.~A.}
\newblock Simulation of multifractal brownian motion.
\newblock {\em Technical report\/} (1998).

\bibitem{cheridito:2003}
{\sc Cheridito, P.}
\newblock Arbitrage in fractional {B}rownian motion models.
\newblock {\em Finance Stoch. 7}, 4 (2003), 533--553.

\bibitem{Coeurjolly:2001}
{\sc Coeurjolly, J.~F.}
\newblock Estimating the parameters of a fractional brownian motion by discrete
  variations of its sample paths.
\newblock {\em Statist. Inf. Stoch. Proc. 4\/} (2001), 199--227.

\bibitem{Coeurjolly:2005}
{\sc Coeurjolly, J.-F.}
\newblock Identification of multifractional {B}rownian motion.
\newblock {\em Bernoulli 11}, 6 (2005), 987--1008.

\bibitem{Cohen:1999}
{\sc Cohen, S.}
\newblock {\em From self-similarity to local self-similarity : the estimation
  problem, Fractal: Theory and Applications in Engineering}.
\newblock Dekking, M., Lévy Véhel, J., Lutton, E., and Tricot, C. (Eds).
  Springer Verlag, 1999.

\bibitem{Cramer:Leadbetter:1967}
{\sc Cram\`{e}r, H., and Leadbetter, M.~R.}
\newblock {\em Stationary and Related Stochastic Processes (Sample Function
  Properties and Their Applications)}.
\newblock 1967.

\bibitem{Fhima:2011}
{\sc Fhima, M.}
\newblock Phd thesis in preparation, 2011.

\bibitem{Kolmogorov:1940}
{\sc Kolmogorov, A.~N.}
\newblock Wienersche spiralen und einige andere interessante kurven im
  hilbertschen raum.
\newblock {\em C. R. (Doklady) Acad. URSS (N.S.) 26\/} (1940), 115--118.

\bibitem{LevyVehel:Peltier:1996}
{\sc Lévy-Véhel, J., and Peltier, R.~F.}
\newblock Multifractional brownian motion : definition and preliminary results.
\newblock {\em Techn. Report RR-2645, INRIA\/} (1996).

\bibitem{Mandelbrot:VanNess:1968}
{\sc Mandelbrot, B., and Van~Ness, J.}
\newblock Fractional brownian motions, fractional noises and applications.
\newblock {\em SIAM 10\/} (1968), 422--437.

\bibitem{Meyer:etal:1999}
{\sc Meyer, Y., Sellan, F., and Taqqu, M.}
\newblock Wavelets, generalized white noise and fractional integration: the
  synthesis of fractional brownian motion.
\newblock {\em J Fourier Anal Appl 5\/} (1999), 465--494.

\bibitem{Nourdin:etal:2010}
{\sc Nourdin, I., Peccati, G., and Podolskij, M.}
\newblock Quantitative breuer-major theorems.
\newblock {\em HAL : hal-00484096, version 2\/} (2010).

\bibitem{Samorodnitsky:Taqqu:1994}
{\sc Samorodnitsky, G., and Taqqu, M.~S.}
\newblock {\em Stable non-Gaussian random processes}.
\newblock Chapman \& Hall, 1994.

\bibitem{Stoncelis:Vaiciulis:2008}
{\sc Stoncelis, M., and Vai\v{c}iulis, M.}
\newblock Numerical approximation of some infinite gaussian series and
  integrals.
\newblock {\em Nonlinear Analysis: Modelling and Control 13\/} (2008),
  397--415.

\bibitem{Surgailis:etal:2008}
{\sc Surgailis, D., Teyssière, G., and Vai\v{c}iulis, M.}
\newblock The increment ratio statistic.
\newblock {\em J. Multivariate Anal 99\/} (2008), 510--541.

\bibitem{Yaglom:1957}
{\sc Yaglom, A.~M.}
\newblock Some classes of random fields in n-dimensional space, related to
  stationary random processes.
\newblock {\em Th. Probab. Appl. 2\/} (1957), 273--320.

\end{thebibliography}


\end{document}